\documentclass[11pt]{article}
\usepackage{amsmath, amssymb, amsthm, amsfonts}
\usepackage{amssymb}
\usepackage{amscd}
\usepackage{verbatim}
\begin{document}
\newcommand{\dx}{\,\mathrm{d}x}
\newcommand{\dy}{\,\mathrm{d}y}
\newcommand{\dz}{\,\mathrm{d}z}
\newcommand{\dt}{\,\mathrm{d}t}
\newcommand{\core}{C_0^{\infty}(\Omega)}
\newcommand{\sob}{W^{1,p}(\Omega)}
\newcommand{\sobloc}{W^{1,p}_{\mathrm{loc}}(\Omega)}
\newcommand{\merhav}{{\mathcal D}^{1,p}}
\newcommand{\be}{\begin{equation}}
\newcommand{\ee}{\end{equation}}
\newcommand{\mysection}[1]{\section{#1}\setcounter{equation}{0}}
\newcommand{\bea}{\begin{eqnarray}}
\newcommand{\eea}{\end{eqnarray}}
\newcommand{\bean}{\begin{eqnarray*}}
\newcommand{\eean}{\end{eqnarray*}}
\newcommand{\thkl}{\rule[-.5mm]{.3mm}{3mm}}
\newcommand{\cw}{\stackrel{\rightharpoonup}{\rightharpoonup}}
\newcommand{\id}{\operatorname{id}}
\newcommand{\supp}{\operatorname{supp}}
\newcommand{\wlim}{\mbox{ w-lim }}
\newcommand{\mymu}{{x_N^{-p_*}}}
\newcommand{\Reals}{{\mathbb R}}
\newcommand{\N}{{\mathbb N}}
\newcommand{\Z}{{\mathbb Z}}
\newcommand{\Q}{{\mathbb Q}}
\newcommand{\CC}{\mathbb{C}}
\newcommand{\Hess}{\mathrm{Hess}}
\newcommand{\Real}{{\mathbb R}}
\newcommand{\R}{{\mathbb R}}
\newcommand{\Rn}{\mathbb{R}^n}
\newcommand{\Zn}{\Z^n}
\newcommand{\Nat}{{\mathbb N}}
\newcommand{\abs}[1]{\lvert#1\rvert}
\newcommand{\Green}[4]{\mbox{$G^{#1}_{#2}(#3,#4)$}}
\newtheorem{theorem}{Theorem}[section]
\newtheorem{corollary}[theorem]{Corollary}
\newtheorem{lemma}[theorem]{Lemma}
\newtheorem{definition}[theorem]{Definition}
\newtheorem{remark}[theorem]{Remark}
\newtheorem{remarks}[theorem]{Remarks}
\newtheorem{proposition}[theorem]{Proposition}
\newtheorem{problem}[theorem]{Problem}
\newtheorem{conjecture}[theorem]{Conjecture}
\newtheorem{question}[theorem]{Question}
\newtheorem{example}[theorem]{Example}
\newtheorem{Thm}[theorem]{Theorem}
\newtheorem{Lem}[theorem]{Lemma}
\newtheorem{Pro}[theorem]{Proposition}
\newtheorem{Def}[theorem]{Definition}
\newtheorem{Exa}[theorem]{Example}
\newtheorem{Exs}[theorem]{Examples}
\newtheorem{Rems}[theorem]{Remarks}
\newtheorem{rem}[theorem]{Remark}
\newtheorem{Cor}[theorem]{Corollary}
\newtheorem{Conj}[theorem]{Conjecture}
\newtheorem{Prob}[theorem]{Problem}
\newtheorem{Ques}[theorem]{Question}
\newcommand{\pf}{\noindent \mbox{{\bf Proof}: }}
\newcommand{\dnorm}[1]{\thkl #1 \thkl\,}

\renewcommand{\theequation}{\thesection.\arabic{equation}}
\catcode`@=11 \@addtoreset{equation}{section} \catcode`@=12
\def\ga{\alpha}     \def\gb{\beta}       \def\gg{\gamma}
\def\gc{\chi}       \def\gd{\delta}      \def\ge{\epsilon}
\def\gth{\theta}                         \def\vge{\varepsilon}
       \def\vgf{\varphi}    \def\gh{\eta}
\def\gi{\iota}      \def\gk{\kappa}      \def\gl{\lambda}
\def\gm{\mu}        \def\gn{\nu}         \def\gp{\pi}
\def\vgp{\varpi}    \def\gr{\rho}        \def\vgr{\varrho}
\def\gs{\sigma}     \def\vgs{\varsigma}  \def\gt{\tau}
\def\gu{\upsilon}   \def\gv{\vartheta}   \def\gw{\omega}
        \def\gy{\psi}        \def\gz{\zeta}
\def\Gg{\Gamma}     \def\Gd{\Delta}      \def\Gf{\Phi}
\def\Gth{Theta}
\def\Gl{\Lambda}    \def\Gs{\Sigma}      \def\Gp{\Pi}
\def\Gw{\Omega}     \def\Gx{\Xi}         \def\Gy{\Psi}


\title{Topics in the theory of positive solutions of second-order elliptic
and parabolic partial differential equations}
\author{Yehuda Pinchover\\
 {\small Department of Mathematics}\\ {\small  Technion - Israel Institute of Technology}\\
 {\small Haifa 32000, Israel}\\
{\small pincho@techunix.technion.ac.il}\\[3mm]
 {\em Dedicated to Barry Simon}\\ {\em on the occasion of his 60th
 birthday}}
\date{}
\maketitle
\begin{abstract}
The purpose of the paper is to review a variety of recent
developments in the theory of positive solutions of general linear
elliptic and parabolic equations of second-order on noncompact
Riemannian manifolds, and to point out a number of their
consequences.
\\[1mm]
\noindent  2000 {\em Mathematics Subject Classification.}
Primary 35J15; Secondary  35B05, 35C15, 35K10.\\[1mm]
 \noindent {\em Keywords.}
Green function, ground state, heat kernel, Liouville theorem,
Martin boundary, positive solution, $p$-Laplacian.
\end{abstract}
 \mysection{Introduction}\label{secint}
Positivity properties of general linear second-order elliptic and
parabolic equations have been extensively studied over the recent
decades (see for example \cite{M98,Pins95} and the references
therein). The purpose of the present paper is to review a variety
of recent developments in the theory of positive solutions of such
equations and to point out a number of their (sometimes
unexpected) consequences. The attention is focused on
generalizations of positivity properties which were studied by
Barry Simon in the special case of Schr\"odinger operators. Still,
the selection of topics in this survey is incomplete, and is
according to the author's working experience and taste. The
reference list is far from being complete and serves only this
expos\'{e}.

The outline of the paper is as follows. In
Section~\ref{secpreliminaries}, we introduce some fundamental
notions that will be studied throughout the paper. In particular,
we bring up the notions of the generalized principal eigenvalue,
criticality and subcriticality of elliptic operators, and the
Martin boundary. Section~\ref{secpert} is devoted to different
types of perturbations and their properties. In
Section~\ref{secindef}, we study the behavior of critical
operators under indefinite perturbations. In
sections~\ref{sechetk} and \ref{secup} we discuss some
relationships between criticality theory and the theory of
nonnegative solutions of the corresponding parabolic equations.
More precisely, in Section~\ref{sechetk} we deal with the large
time behavior of the heat kernel, while in Section~\ref{secup} we
discuss sufficient conditions for the nonuniqueness of the
positive Cauchy problem, and study intrinsic ultracontractivity.

In Section~\ref{seceigen}, we study the asymptotic behavior at
infinity of eigenfunctions of Schr\"odinger operators. The
phenomenon known in the mathematical physics literature as
`localization of binding', and the properties of the shuttle
operator are discussed in sections~\ref{seclocalization} and
\ref{secshttle}, respectively. The exact asymptotics of the
positive minimal Green function, and the explicit Martin integral
representation theorem for positive solutions of general
$\mathbb{Z}^d$-periodic elliptic operators on $\mathbb{R}^d$ are
reviewed in Section~\ref{secperiod}. We devote
Section~\ref{seccritliouville} to some relationships between
criticality theory and Liouville theorems. In particular, we
reveal that an old open problem of B.~Simon
(Problem~\ref{problem}) is completely solved (see
Theorem~\ref{thmDKS}). In Section~\ref{secliouville} we study
polynomially growing solutions of $\mathbb{Z}^d$-periodic
equations on $\mathbb{R}^d$. We conclude the paper in
Section~\ref{secplap} with criticality theory for the
$p$-Laplacian with a potential term.
\mysection{Principal eigenvalue, minimal growth  and
classification}\label{secpreliminaries}
  Consider a noncompact, connected, smooth Riemannian manifold $X$ of dimension $d$.
For any subdomain $\Omega\subseteq X$, we write $D\Subset \Omega$
if $\overline{D}$ is a compact subset of $\Omega$. The ball of
radius $r>0$ and center at $x_0$ is denoted by $B(x_0,r)$. Let
$f,g \in C(\Omega)$, we use the notation $f\asymp g$ on
$D\subseteq\Omega$ if there exists a positive constant $C$ such
that
$$C^{-1}g(x)\leq f(x) \leq Cg(x) \qquad \mbox{ for all } x\in D.$$
By ${\bf 1}$, we denote the constant function taking at any point
the value $1$.

We associate to any subdomain $\Omega\subseteq X$ {\em an
exhaustion of $\Omega$}, i.e. a sequence of smooth, relatively
compact domains $\{\Omega_{j}\}_{j=1}^{\infty}$ such that
$\Omega_1\neq \emptyset$, $\overline{\Omega}_{j}\subset
\Omega_{j+1}$ and $\cup_{j=1}^{\infty}\Omega_{j}=\Omega$. For
every $j\geq 1$, we denote $\Omega_{j}^*=\Omega\setminus
\overline{\Omega_j}$. We say that a function $f\in C(\Omega)$ {\em
vanishes at infinity of} $\Omega$ if for  every $\varepsilon >0$
there exists $N\in \Nat$ such that $|f(x)|<\varepsilon$ for all $x
\in \Omega^*_{N}$.

We associate to any such exhaustion
$\{\Omega_{j}\}_{j=1}^{\infty}$ a sequence
$\{\chi_{j}(x)\}_{j=1}^{\infty}$ of smooth cutoff functions in
$\Omega$ such that $\chi_j(x)\equiv 1$ in $\Omega_j$,
$\chi_j(x)\equiv 0$ in $\Omega \setminus \Omega_{j+1}$, and $0\leq
\chi_j(x)\leq 1$ in $\Omega$. Let $0<\ga \leq 1$. For $W\in
C^\ga(\Omega)$,  we denote $W_j(x)=\chi_j(x)W(x)$ and
$W^*_j(x)=W(x)-W_j(x)$.

\indent We consider a linear, second-order, elliptic operator $P$
defined in a subdomain $\Omega\subset X$. Here $P$ is an operator
with {\em real} H\"{o}lder continuous coefficients which in any
coordinate system $(U;x_{1},\ldots,x_{d})$ has the form
\begin{equation} \label{P} P(x,\partial_{x})=-\sum_{i,j=1}^{d}
a_{ij}(x)\partial_{i}\partial_{j} + \sum_{i=1}^{d}
b_{i}(x)\partial_{i}+c(x), \end{equation} where
$\partial_{i}=\partial/\partial x_{i}$. We assume that for each
$x\in \Omega$ the real quadratic form $\sum_{i,j=1}^{d}
a_{ij}(x)\xi_{i}\xi_{j}$ is positive definite on $\mathbb{R} ^d$.

We denote the cone of all positive (classical) solutions of the
elliptic equation $Pu=0$ in $\Omega$ by $\mathcal{C}_{P}(\Omega)$.
We fix a reference point  $x_0 \in \Omega _1$. From time to time,
we consider the convex set
$$\mathcal{K}_{P}(\Omega):=\{u\in \mathcal{C}_{P}(\Omega)\,|\,u(x_0)=1\}$$
of all {\em normalized} positive solutions. In case that the
coefficients of $P$ are smooth enough, we denote by $P^*$ the
formal adjoint of $P$.

\begin{definition}{\em For a (real valued) function $V\in C^\alpha(\Omega)$, let
$$\lambda_0(P,\Omega,V)
:= \sup\{\lambda \in \mathbb{R} \mid \mathcal{C}_{P-\lambda V}(
\Omega)\neq \emptyset\}$$ be the {\em generalized principal
eigenvalue} of the operator $P$ with respect to the (indefinite)
weight $V$ in $\Omega$.  We also denote
 $$\lambda_\infty(P,\Omega, V):=
 \sup_{K\Subset\Omega}\lambda_0(P,\Omega\setminus K,V).$$
For a fixed $P$ and $\Omega$, and $V=\mathbf{1}$, we simply write
$\lambda_0:=\lambda_0(P,\Omega,\mathbf{1})$ and
$\lambda_\infty:=\lambda_\infty(P,\Omega,\mathbf{1})$.
 }\end{definition}
\begin{definition}{\em
Let  $P$ be an elliptic operator of the form (\ref{P}) which is
defined on a smooth domain $D\Subset X$. we say that the {\em
generalized maximum principle} for the operator $P$ holds in $D$
if for any $u\in C^2(D)\cap C(\overline{D})$, the inequalities
$Pu\geq 0$ in $D$ and $u\geq 0$ on $\partial D$ imply that $u\geq
0$ in $D$.
 }\end{definition}
It is well known that $\lambda_0(P,\Omega,\mathbf{1})\geq 0$ if
and only if the generalized maximum principle for the operator $P$
holds true in any smooth subdomain $D\Subset\Omega$.

 The following theorem is known as the {\em Allegretto-Piepenbrink
theory},  it relates $\lambda_0$ and $\lambda_\infty$, in the
symmetric case,  with fundamental spectral quantities (see for
example \cite{Agmon82,Cycon,S82} and the references therein).
\begin{theorem}\label{APthm} Suppose that $P$ is symmetric on
$C_0^\infty(\Omega)$, and that $\lambda_0>-\infty$. Then
$\lambda_0$ (resp. $\lambda_\infty$) equals to the infimum of the
spectrum (resp. essential spectrum) of the Friedrich's extension
of $P$.
 \end{theorem}
Therefore, in the selfadjoint case, $\lambda_0$ can be
characterized via the classical Rayleigh-Ritz variational formula.
In the general case, a variational principle for $\lambda_0$ is
given by the Donsker-Varadhan variational formula (which is a
generalization of the Rayleigh-Ritz formula) and by some other
variational formulas (see for example \cite{NP,Pins95}).

\begin{Def} \label{defminimalg}
{\rm Let $P$ be an elliptic operator defined in a domain $\Omega
\subseteq X$. A function $u$ is said to be a {\em positive
solution of the operator $P$ of minimal growth in a neighborhood
of infinity in} $\Omega$ if $u\in \mathcal{C}_P(\Omega _j^*)$ for
some $j\geq 1$, and for any  $l> j$, and $v\in C(\Omega_l^*\cup
\partial \Omega _l)\cap \mathcal{C}_P(\Omega _l^*)$, if $u\le v$
on $\partial \Omega _l$, then $u\le v$ on $\Omega _l^*$.}
\end{Def}

\begin{theorem}[\cite{Agmon82}]\label{thmmingr2}
Suppose that $\mathcal{C}_{P}(\Omega)\neq \emptyset$. Then for any
$x_0\in \Omega$ the equation $Pu=0$ has (up to a multiple
constant) a unique positive solution $v$ in
$\Omega\setminus\{x_0\}$ of minimal growth in a neighborhood of
infinity in $\Omega$.
\end{theorem}
By the well known theorem on the removability of isolated
singularity \cite{GiS}, we have:
\begin{Def} {\em
Suppose that $\mathcal{C}_{P}(\Omega)\neq \emptyset$. If the
solution $v$ of Theorem~\ref{thmmingr2} has a nonremovable
singularity at $x_0$, then $P$ is said to be a {\em subcritical
operator} in $\Omega$. If $v$ can be (uniquely) continued to a
positive solution $\tilde{v}$ of the equation $Pu=0$ in $\Omega$,
then $P$ is said to be a {\em critical operator} in $\Omega$, and
the positive global solution  $\tilde{v}$   is called a {\em
ground state} of the equation $Pu=0$ in $\Omega$. The operator $P$
is said to be {\em supercritical} in $\Omega$ if
$\mathcal{C}_{P}(\Omega)=\emptyset$.
 } \end{Def}

\begin{remarks}{\em
1. In \cite{S80}, B.~Simon coined the terms
`(sub)-(super)-critical operators' for Schr\"odinger operators
with short-range potentials which are defined on $\mathbb{R}^d$,
where $d\geq 3$. The definition given in \cite{S80} is in terms of
the exact (and particular) large time behavior of the heat kernel
of such operators (see \cite[p.~71]{S81} for the root of this
terminology). In \cite{M86}, M.~Murata generalized the above
classification for Schr\"odinger operators which are defined in
any  subdomain of $\mathbb{R}^d$, $d\geq 1$. The definition of
subcriticality given here is due to \cite{P88}.

2. The notions of minimal growth and ground state were introduced
by S.~Agmon in \cite{Agmon82}.

3. For modified and stronger notions of subcriticality see
\cite{DS91,P88}.

 }\end{remarks}

\noindent {\em Outline of the proof of Theorem~\ref{thmmingr2}.}
Assume that $\mathcal{C}_{P}(\Omega)\neq \emptyset$ and fix
$x_0\in \Omega$. Then for every $j\geq 1$, the {\em Dirichlet
Green function} $\Green{\Omega_{j}}{P}{x}{y}$ for the operator $P$
exists in $\Omega_j$. It is the integral kernel such that for any
$f\in\core$, the function $u_j(x):=\int_{\Omega_j}
G_P^{\Omega_j}(x,y)f(y)\, \mathrm{d}y$ solves the Dirichlet
boundary value problem
$$
 Pu=f \quad \mbox{ in } \Omega_j,\qquad
  u=0 \quad \mbox{ on  } \partial \Omega_j.
 $$
It follows that $\Green{\Omega_{j}}{P}{\cdot}{x_0}\in
\mathcal{C}_{P}(\Omega_j\setminus\{x_0\})$. By the generalized
maximum principle,
$\{\Green{\Omega_{j}}{P}{x}{x_0}\}_{j=1}^{\infty}$ is an
increasing sequence which,  by the Harnack inequality, converges
uniformly in any compact  subdomain of $\Omega\setminus\{x_0\}$
either to $\Green{\Omega}{P}{x}{x_0}$, the  positive  {\em minimal
Green function} of $P$ in $\Omega$ with a pole at $x_0$ (and in
this case $P$ is subcritical in $\Omega$)  or to infinity.

In the latter case, fix $x_1\in \Omega$, such that $x_1\neq x_0$.
It follows that the sequence
$\Green{\Omega_j}{P}{\cdot}{x_0}/\Green{\Omega_j}{P}{x_1}{x_0}$
converges uniformly in any compact  subdomain of
$\Omega\setminus\{x_0\}$ to a ground state of the equation $Pu=0$
in $\Omega$, and in this case $P$ is critical in $\Omega$.\qed

\begin{corollary} \label{cor27}
 (i) If $P$ is subcritical in $\Omega$, then for each $y\in \Omega$ the
Green function $\Green{\Omega}{P}{\cdot}{y}$ with a pole at $y$
exists, and is a positive solution of the equation $Pu=0$ of
minimal growth in a neighborhood of infinity in $\Omega$.
Moreover, $P$ is subcritical in $\Omega$ if and only if the
equation $Pu=0$ in $\Omega$ admits a positive supersolution which
is not a solution.

(ii) The operator $P$ is critical in $\Omega$ if and only if the
equation $Pu=0$ in $\Omega$ admits (up to a multiplicative
constant) a unique positive supersolution. In particular, $\dim
\mathcal{C}_{P}(\Omega)=1$.

(iii) Suppose that $P$ is symmetric on $C_0^\infty(\Omega)$  with
respect to a smooth positive density $V$, and let $\tilde{P}$ be
the (Dirichlet) selfadjoint realization of $P$ on
$L^2(\Omega,V(x)dx)$. Assume that
$\lambda\in\sigma_\mathrm{point}(\tilde{P})$ admits a nonnegative
eigenfunction $\varphi$, then $\lambda\!=\!\lambda_0$ and
$P\!-\!\lambda_0 V$ is critical in $\Omega$ (see for example
\cite{M86}).

 (iv) The operator $P$ is critical (resp. subcritical) in
$\Omega$ if and only if $P^*$ is critical (resp. subcritical) in
$\Omega$.
\end{corollary}
As was mentioned, (sub)criticality is related to the large time
behavior of the heat kernel. Indeed, (sub)criticality can be also
defined in terms of the corresponding parabolic equation. Suppose
that $\lambda_0\geq 0$. For every $j\geq 1$, consider the
Dirichlet heat kernel $k_P^{\Omega_j}(x,y,t)$ of the parabolic
operator $L:=\partial_t+P$ on $\Omega_j\times (0,\infty)$. So, for
any $f\in\core$, the function $u_j(x,t)=\int_{\Omega_j}
k_P^{\Omega_j}(x,y,t)f(y)\, \mathrm{d}y$ solves the
initial-Dirichlet boundary value problem
 $$
 Lu=0 \; \mbox{ in } \Omega_j\times (0,\infty),\quad
  u=0 \; \mbox{ on  } \partial \Omega_j\times (0,\infty),\quad
  u=f  \; \mbox{ on  }  \Omega_j\times \{0\}.
$$
   By the (parabolic) generalized maximum principle, $\{k_P^{\Omega_j}(x,y,t)\}_{j=1}^{\infty}$ is
an increasing sequence which converges to $k_P^{\Omega}(x,y,t)$,
the {\em minimal heat kernel} of the parabolic operator $L$ in
$\Omega$.
\begin{lemma}\label{lemheatcrit}
Suppose that $\lambda_0\geq 0$. Let $x,y\in \Omega$, $x\neq y$.
Then
$$\int_0^\infty k_P^{\Omega}(x,y,t)\,\mathrm{d}t<\infty \qquad
\mbox{(resp. $\int_0^\infty
k_P^{\Omega}(x,y,t)\,\mathrm{d}t=\infty$),}$$ if and only if  $P$
is a subcritical (resp. critical) operator in $\Omega$. Moreover,
if $P$ is subcritical operator in $\Omega$, then
\begin{equation}\label{heatgreen}G_P^{\Omega}(x,y)=\int_0^\infty
k_P^{\Omega}(x,y,t)\,\mathrm{d}t.
 \end{equation}
\end{lemma}

For the proof of Lemma~\ref{lemheatcrit} see for example
\cite{Pins95}. Note that if $\lambda<\lambda_0$, then the operator
$P-\lambda$ is subcritical in $\Omega$, and that for $\lambda\leq
\lambda_0$, the heat kernel $k_{P-\lambda}^\Omega(x,y,t)$ of the
operator $P-\lambda$ is equal to $e^{\lambda t}k_P^\Omega(x,y,t)$.

\vskip 3mm

Subcriticality (criticality) can be defined also through a
probabilistic approach. If the zero-order coefficient $c$ of the
operator $P$ is equal to zero in $\Omega$, then $P$ is called a
{\em diffusion operator}. In this case, $P\mathbf{1}=0$, and
therefore, $P$ is not supercritical in $\Omega$. Moreover, for
such an operator $P$, one can associate a diffusion process
corresponding to a solution of the generalized martingale problem
for $P$ in $\Omega$. This diffusion process is either {\em
transient} or {\em recurrent} in $\Omega$. It turns out that a
diffusion operator $P$ is subcritical in $\Omega$ if and only if
the associated diffusion process is transient in $\Omega$ (for
more details see \cite{Pins95}). A Riemannian manifold $X$ is
called {\em parabolic} (resp. {\em non-parabolic}) if the Brwonian
motion, the diffusion process with respect to the Laplace-Beltrami
operator on $X$, is recurrent (resp. transient) \cite{G99}.

Suppose now that $P$ is of the form \eqref{P}, and $P$ is  not
supercritical in $\Omega$. Let $\varphi\in
\mathcal{C}_{P}(\Omega)$. Then the operator $P^\varphi$ acting on
functions $u$ by
$$P^\varphi u:=\frac{1}{\varphi}P(\varphi u)$$
is a diffusion operator, and $$
k_{P^\varphi}^M(x,y,t)=\frac{1}{\varphi(x)}k_{P}^M(x,y,t)\varphi(y).$$
Therefore, $P$ is subcritical in $\Omega$ if and only if
$P^\varphi$ is transient in $\Omega$.

\vskip 3mm

We have the following general convexity results.
\begin{theorem}[\cite{P90}]\label{thmconv}
(i) Let $V\in C^\ga(\Omega)$, $V\neq 0$ and set
 \bea S_+ & = &
S_+(P,\Omega,V)  =
\{\lambda\in \mathbb{R}\;|\, P-\lambda V \mbox{ is subcritical in } \Omega\},\\
S_0 & = & S_0(P,\Omega,V)  = \{\lambda\in \mathbb{R}\;|\,
P-\lambda V \mbox{ is critical in } \Omega\}.
 \eea
  Then  $S:=S_+\cup S_0\subseteq
\mathbb{R}$ is a closed interval and $S_0\subset \partial S$.
Moreover, if $S\neq \emptyset$, then $S$ is bounded if and only if
$V$ changes its sign in $\Omega$. \vskip 2mm

(ii) Let $W, V\in C^\ga(\Omega)$, then the function
$\lambda_0(\mu):=\lambda_0(P-\mu W,\Omega,V)$ is a concave
function on the interval $\{\mu\in \mathbb{R}\mid
|\lambda_0(\mu)|<\infty\}$.
\end{theorem}
\begin{proof} For $0\leq s\leq 1$, and $V_0,V_1\in C^\ga(\Omega)$,
let $$P_s:= P+sV_1+(1-s)V_0.$$  Assume that $u_j$ are positive
supersolutions of the  equations $P_j u\geq 0$ in $\Omega$, where
$j=0,1$. It can be verified that for $0< s<1$, the function
$$u_s(x):=\left[u_0(x)\right]^{1-s}\left[u_1(x)\right]^s$$
is a positive supersolution of the equation $P_su=0$ in $\Omega$.
Moreover, for any $0< s<1$, $u_s\in \mathcal{C}_{P_s}(\Omega)$ if
and only if $V_0=V_1$, and $u_0,u_1\in \mathcal{C}_{P_0}(\Omega)$
are linearly dependent. The lemma follows easily from this
observation.
\end{proof}
\begin{corollary}[\cite{P90} and \cite{S82}]\label{cor3G}
Suppose that $P_s:= P+sV_1+(1-s)V_0$ is subcritical in $\Omega$
for $s=0,1$. Then for $0\leq s\leq 1$ we have
\begin{equation}\label{eq3G}
G_{P_s}^\Omega(x,y)\leq
\left[G_{P_0}^\Omega(x,y)\right]^{1-s}\left[G_{P_1}^\Omega(x,y)\right]^{s}.
\end{equation}
\end{corollary}
\begin{remark}{\em
The dependence of $\lambda_0$ on the higher order coefficients of
$P$ is more involved. In \cite{BNV} it was proved that in the
class of uniformly elliptic operators with bounded coefficients
which are defined on a bounded domain in $\mathbb{R}^d$,
$\lambda_0$ is locally Lipschitz continuous as a function of the
first-order coefficients of the operator $P$. A.~Ancona
\cite{An97} proved that under some assumptions, $\lambda_0$ is
Lipschitz continuous with respect to a metric
$\mathrm{dist}(P_1,P_2)$ measuring the distance between two
elliptic operators $P_1$ and $P_2$ in a certain class. Ancona's
metric depends on the difference between \emph{all} the
coefficients of the operators $P_1$ and $P_2$. }
\end{remark}

\vskip 3mm

If $P$ is subcritical in $\Omega$, then $\mathcal{C}_{P}(\Omega)$
is in general not a one-dimensional cone. Nevertheless, one can
construct the {\em Martin compactification} $\Omega_P^M$ of
$\Omega$ with respect to the operator $P$ (with a base point
$x_0$), and obtain an integral representation of any solution in
$\mathcal{C}_{P}(\Omega)$.  More precisely, the {\em Martin
compactification} is the compactification of $\Omega$ such that
the function
 $$K_P^\Omega(x,y):=\frac{G_P^\Omega(x,y)}{G_P^\Omega(x_0,y)}
\qquad \mbox{ on } \Omega\times \Omega \setminus\{(x_0,x_0)\}$$
has a continuous extension $K_P^\Omega(x,\eta)$
 to
$\Omega\times (\Omega_P^M\setminus \{x_0\})$, and such that the
set of functions  $\{K_P^\Omega(\cdot,\eta)\}_{\eta\in
\Omega_P^M}$ separates the points of $\Omega_P^M$.
 The boundary of $\Omega_P^M$ is
denoted by $\partial_P^M \Omega$ and is called the {\em Martin
boundary} of $\Omega$ with respect to the operator $P$. For each
$\xi \in
\partial _P^M\Omega$, the function
$K_P^\Omega(\cdot,\xi)$ is called the {\em Martin function of the
pair $(P,\Omega)$ with a pole at $\xi$}. Note that for $\xi \in
\partial _P^M\Omega$, we have
$K_P^\Omega(\cdot,\xi)\in \mathcal{K}_{P}(\Omega)$.  The set
$\partial _{m,P}^M\Omega$ of all $\xi \in
\partial _P^M\Omega$ such that $K_P^\Omega(\cdot,\xi)$ is an extreme point of
the convex set $\mathcal{K}_{P}(\Omega)$ is called the {\em
minimal Martin boundary}  (for more details see \cite[and the
references therein]{M86,M98,Pins95,Taylor}).

The Martin representation theorem asserts that for any $u\in
\mathcal{K}_{P}(\Omega)$ there exists a unique probability measure
$\mu$ on $\partial_P^M \Omega$ which is supported on $\partial
_{m,P}^M\Omega$ such that
$$u(x)=\int_{\partial_P^M \Omega}K_P^\Omega(x,\xi)\,\mathrm{d}\mu(\xi).$$

 There has been a great deal of work on explicit description of the Martin
compactification and representation in many concrete examples (see
for example \cite[and the references
therein]{LP,M98,MT,Pins95,Taylor}).

We present below two elementary examples of Martin
compactifications. In Section~\ref{secperiod} we discuss a recent
result on the Martin compactification of a general periodic
operator on $\mathbb{R}^d$.
\begin{example}\label{exsmoothbdd} {\em
Let $\Omega$ be a smooth bounded domain in $\mathbb{R}^d$, and
assume that the coefficients of $P$ are (up to the boundary)
smooth. Then $\partial_P^M \Omega$ is homeomorphic to $\partial
\Omega$, the euclidian boundary of $\Omega$, and for any $y\in
\partial \Omega$,
\begin{equation}\label{eqmartincalss}
 K_P^\Omega(x,y):=\frac{\partial_\nu
G_P^\Omega (x,y)} {\partial_\nu G_P^\Omega(x_0,y)}\,,
\end{equation}
where $\partial_\nu$ denotes the inner normal derivative with
respect to the second variable. Note that $\partial_\nu
G_P^\Omega(\cdot,y)$ is the Poisson kernel at $y\in
\partial \Omega$. }
\end{example}
\begin{example}\label{exlaplace} {\em
 Consider the equation $H_\lambda u:=(-\Delta+\lambda)u=0$ in
$\mathbb{R}^d$. Then $\mathcal{C}_{H_\lambda}(\mathbb{R}^d)\neq
\emptyset$ if and only if $\lambda\geq 0$. It is well known that
$H_0=-\Delta$ is critical on in $\mathbb{R}^d$ if and only if
$d\leq 2$. Moreover,
$$G_{H_\lambda}^{\mathbb{R}^d}(x,y)=
  \begin{cases}
    \dfrac{\Gamma(\nu)|x-y|^{2-d}}{4\pi^{d/2}} & \lambda=0, \text{ and } d\geq 3, \\[3mm]
     (2\pi)^{-d/2}\left(\dfrac{\sqrt{\lambda}}{|x-y|}\right)^\nu K_\nu(\sqrt{\lambda}|x-y|) & \lambda>0,
  \end{cases}$$
where $\nu=(d-2)/2$, and $K_\nu$ is the modified Bessel function
of order $\nu$.

Clearly,
$$\lim_{|y|\to \infty}
\dfrac{G_{-\Delta}^{\mathbb{R}^d}(x,y)}{G_{-\Delta}^{\mathbb{R}^d}(0,y)}=1.$$
Therefore, the Martin compactification of $\mathbb{R}^d$ with
respect to the Laplacian is the one-point compactification of
$\mathbb{R}^d$, and we obtained the positive Liouville theorem:
$\mathcal{K}_{-\Delta}(\mathbb{R}^d)=\{\mathbf{1}\}$.

Suppose now that $\lambda>0$. Then for any $\xi\in S^{d-1}$,
$$\lim_{\frac{y}{|y|}\to \xi,\; |y|\to \infty}\;
\dfrac{G_{H_\lambda}^{\mathbb{R}^d}(x,y)}{G_{H_\lambda}^{\mathbb{R}^d}(0,y)}=e^{\sqrt{\lambda}\,\xi\cdot
x},$$ and therefore, the Martin boundary of $\mathbb{R}^d$ with
respect to $H_\lambda$ is the sphere at infinity.
 Clearly, all
Martin functions are minimal. Furthermore, $u\in
\mathcal{C}_{H_\lambda}(\mathbb{R}^d)$ if and only if there exists
a positive finite measure $\gm$ on $S^{d-1}$ such that
\[
u(x)=\int_{S^{d-1}} e^{\sqrt{\lambda}\,\xi\cdot
x}\,\mathrm{d}\gm(\xi ).
\]
 }
\end{example}
\begin{remark}{\em We would like to point out that criticality
theory and Martin boundary theory are also valid  for the class of
weak solutions of elliptic equations in divergence form as well as
for the class of strong solutions of strongly elliptic equations
with locally bounded coefficients. For the sake of clarity, we
prefer to concentrate on the class of classical solutions.}
\end{remark}
\mysection{Perturbations}\label{secpert}
An operator $P$ is critical in $\Omega$ if and only if any
positive supersolution of the equation $Pu=0$ in $\Omega$ is a
solution (Corollary \ref{cor27}). Therefore, if $P$ is critical in
$\Omega$ and $V\in C^\ga(\Omega)$ is a nonzero, nonnegative
function, then for any $\lambda>0$ the operator $P+\lambda V$ is
subcritical and $P-\lambda V$ is supercritical in $\Omega$. On the
other hand, it can be shown that subcriticality is a stable
property in the following sense: if $P$ is subcritical in $\Omega$
and $V\in C^\ga(\Omega)$ has a {\em compact support}, then there
exists $\epsilon>0$ such that $P-\lambda V$ is subcritical for all
$|\lambda|< \ge$, and the Martin compactifications $\Omega_{P}^M$
and $\Omega_{P-\lambda V}^M$ are homeomorphic for all $|\lambda|<
\ge$ (for a more general result see Theorem~\ref{thmssp}).
Therefore, a perturbation by a compactly supported potential (at
least with a definite sign) is well understood.

In this section, we introduce and study a few general notions of
perturbations related to positive solutions of an operator $P$ of
the form \eqref{P} by a (real valued) potential $V$. In
particular, we discuss the behavior of the generalized principal
eigenvalue, (sub)criticality, the Green function, and the Martin
boundary under such perturbations. Further aspects of perturbation
theory will be discussed in the following sections.

One facet of this study is the equivalence (or comparability) of
the corresponding Green functions.
\begin{Def} \label{equiGdef}
{\em Let $P_{j},\;j=1,2$, be two subcritical operators in
$\Omega$. We say that the Green functions $G^\Omega_{P_{1}}$ and
$G^\Omega_{P_{2}}$ are {\em equivalent} (resp. {\em
semi-equivalent}) if $G_{P_1}^\Omega \asymp G_{P_2}^\Omega$  on
$\Omega \times \Omega \setminus \{(x,x)\,|\, x \in \Omega\}$
(resp. $G_{P_1}^\Omega (\cdot,y_0)\asymp
G_{P_2}^\Omega(\cdot,y_0)$ on $\Omega \setminus \{ y_0\}$  for
some fixed $y_0\in \Omega$).
 }\end{Def}
\begin{lemma}[\cite{P88}]\label{lemgreeneq} Suppose that
the Green functions $G^\Omega_{P_{1}}$ and $G^\Omega_{P_{2}}$ are
equivalent. Then there exists a homeomorphism $\Phi:\partial
_{m,P_1}^M\Omega\to
\partial _{m,P_2}^M\Omega$ such that for each minimal point $\xi \in
\partial _{m,P_1}^M\Omega$, we have $K_{P_1}^\Omega(\cdot,\xi)
\asymp K_{P_2}^\Omega(\cdot,\Phi(\xi))$ on $\Omega$. Moreover, the
cones $\mathcal{C}_{P_1}(\Omega)$ and $\mathcal{C}_{P_2}(\Omega)$
are homeomorphic.
\end{lemma}

\begin{remarks}{\em
1. It is not known whether the equivalence of $G^\Omega_{P_{1}}$
and $G^\Omega_{P_{2}}$ implies that the cones
$\mathcal{C}_{P_1}(\Omega)$ and $\mathcal{C}_{P_2}(\Omega)$ are
{\em affine} homeomorphic.

2. Many papers deal with sufficient conditions, in terms of
proximity near infinity in $\Omega$ between two given subcritical
operators $P_1$ and $P_2$, which imply that $G^\Omega_{P_{1}}$ and
$G^\Omega_{P_{2}}$ are equivalent, or even that the cones
$\mathcal{C}_{P_1}(\Omega)$ and $\mathcal{C}_{P_2}(\Omega)$ are
affine homeomorphic, see Theorem~\ref{thmssp} and  \cite[and the
references therein]{Ai,An97,M86,M97,P88,P89,Sem}.
 }\end{remarks}
We use the notation
$$E_+=E_+(V,P,\Omega):=\left\{\lambda\in\mathbb{R}\,|\,G^{\Omega}_{P-\lambda V} \mbox{ and }
G^{\Omega}_{P} \mbox{ are equivalent } \right \},$$
$$sE_+=sE_+ (V,P,\Omega):=\left\{\lambda\in\mathbb{R}\,|\,G^{\Omega}_{P-\lambda V}
\mbox{ and } G^{\Omega}_{P} \mbox{ are semi-equivalent }\right
\}.$$

 The following notion was
introduced in \cite{P89} and is closely related to the stability
of $\mathcal{C}_P(\Omega)$ under perturbation by a potential $V$.
\begin{Def} \label{spertdef}
{\em   Let $P$ be a subcritical operator in $\Omega$, and let
$V\in C^{\alpha}(\Omega)$. We say that $V$ is a {\em small
perturbation}
 of $P$ in $\Omega$ if
\be \label{sperteq} \lim_{j\rightarrow \infty}\left\{\sup_{x,y\in
\Omega_{j}^*}
\int_{\Omega_{j}^*}\frac{\Green{\Omega}{P}{x}{z}|V(z)|
\Green{\Omega}{P}{z}{y}}{\Green{\Omega}{P}{x}{y}}\dz\right\}=0.
\ee }
\end{Def}
\noindent The following notions of perturbations were introduced
by M.~Murata \cite{M97}.
\begin{Def} \label{semispertdef}

{\em Let $P$ be a subcritical operator in $\Omega$, and let $V\in
C^{\alpha}(\Omega)$.

{\em (i)} We say that $V$ is a {\em  semismall perturbation} of
$P$ in $\Omega$ if \be \label{semisperteq} \lim_{j\rightarrow
\infty}\left\{\sup_{y\in \Omega_{j}^*} \int_{\Omega_{j}^*}
\frac{\Green{\Omega}{P}{x_0}{z}|V(z)|\Green{\Omega}{P}{z}{y}}
{\Green{\Omega}{P}{x_0}{y}}\dz\right\}=0. \ee

{\em (ii)} We say that $V$ is a {\em  $G$-bounded perturbation}
(resp. {\em $G$-semibounded perturbation}) of $P$ in $\Omega$ if
there exists a positive constant $C$ such that \be \label{bperteq}
\int_{\Omega}\frac{\Green{\Omega}{P}{x}{z}|V(z)|\Green{\Omega}{P}{z}{y}}
{\Green{\Omega}{P}{x}{y}}\dz\leq C \ee for all $x,y \in \Omega$
(resp. for some fixed $x\in \Omega$  and all $y \in \Omega
\setminus \{x\}$).

{\em (iii)} We say that $V$ is an {\em  $H$-bounded perturbation}
(resp. {\em $H$-semibounded perturbation}) of $P$ in $\Omega$ if
there exists a positive constant $C$ such that \be
\label{Hbperteq}
\int_{\Omega}\frac{\Green{\Omega}{P}{x}{z}|V(z)|u(z)}
{u(x)}\dz\leq C \ee for all $x\in \Omega$ (resp. for some fixed
$x\in \Omega$) and all $u \in \mathcal{C}_P(\Omega)$.

{\em (iv)} We say that $V$ is an {\em  $H$-integrable
perturbation} of $P$ in $\Omega$ if \be \label{Iperteq}
\int_{\Omega}\Green{\Omega}{P}{x}{z}|V(z)|u(z)\dz < \infty \ee for
all $x\in \Omega$ and all $u \in \mathcal{C}_P(\Omega)$.
}\end{Def}
\begin{theorem}[\cite{M97,P89,P90}]\label{thmssp}
 Suppose that $P$ is subcritical in $\Omega$.
 Assume that $V$ is a small (resp. semismall) perturbation of
$P^*$ in $\Omega$. Then $E_+=S_+$ (resp. $sE_+=S_+$), and
$\partial S=S_0$. In particular, $S_+$ is an open interval.

Suppose that $V$ is a semismall perturbation of $P^*$ in $\Omega$,
and $\gl\in S_0$. Let $\varphi_0$ be the corresponding ground
state. Then $\varphi_0\asymp \Green{\Omega}{P}{\cdot}{x_0}$ in
$\Omega_1^*$.

Suppose that $V$ is a semismall perturbation of $P^*$ in $\Omega$,
and $\lambda\in S_+$. Then the mapping \be \label{sperteq7}
\Psi(u):=u(x)+\lambda\int_{\Omega}\Green{\Omega}{P-\lambda V
}{x}{z}V(z) u(z)\dz \ee is an affine homeomorphism of
$\mathcal{C}_{P}(\Omega)$ onto $\mathcal{C}_{P-\lambda
V}(\Omega)$, which induces a homeomorphism between the
corresponding Martin boundaries. Moreover, in the  small
perturbation case, we have $\Psi(u)\asymp u$ in $\Omega$ for all
$u\in \mathcal{C}_{P}(\Omega)$.
\end{theorem}

\begin{Rems}\label{remspert}{\em
1. Small perturbations are semismall \cite{M97}, $G$-(resp. $H$-)
bounded perturbations are $G$- (resp. $H$-) semibounded, and
$H$-semibounded perturbations are $H$-integrable. On the other
hand, if $V$ is $H$-integrable and
$\dim\mathcal{C}_P(\Omega)<\infty$, then $V$ is $H$-semibounded
\cite{M97,P88}.

There are potentials which are $H$-semibounded perturbations but
are neither $H$-bounded nor $G$-semibounded. We do not know of any
example of a semismall (resp. $G$-semibounded) perturbation which
is not a small (resp. $G$-bounded) perturbation. We are also not
aware of any example of a $H$-bounded (resp. $H$-integrable)
perturbation which is not $G$-bounded (resp. $H$-semibounded)
\cite{P99}.

2.   Any small (resp. semismall) perturbation  is $G$-bounded
(resp. $G$-semibounded), and any  $G$-(resp. semi) bounded
perturbation is $H$-(resp. semi) bounded perturbation.

3. If $V$ is a $G$-bounded (resp. $G$-semibounded) perturbation of
$P$ (resp. $P^*$) in $\Omega$, then $G^{\Omega}_{P}$ and
$G^{\Omega}_{P-\lambda V}$ are equivalent (resp. semi-equivalent)
provided that $|\lambda |$ is small enough \cite{M97, P88, P89}.
On the other hand, if $G^{\Omega}_{P}$ and $G^{\Omega}_{P+V}$ are
equivalent (resp. semi-equivalent) and $V$ {\em has a definite
sign}, then $V$ is a $G$-bounded (resp. $G$-semibounded)
perturbation of $P$ (resp. $P^*$) in $\Omega$. In this case, by
\eqref{eq3G}, the set $E_+$ (resp. $sE_+$) is an open half line
which is contained in $S_+$ \cite[Corollary 3.6]{P90}. There are
sign-definite $G$-bounded (resp. $G$-semibounded) perturbations
such that $E_+\subsetneqq S_+$ (resp. $sE_+\subsetneqq S_+$)
\cite[Example~8.6]{P99}, \cite[Theorem~6.5]{M98}.

Note that, if $V$ is a $G$-(resp. semi-) bounded perturbation of
$P$ (resp. $P^*$) in $\Omega$ and $\Theta\in C^{\ga}(\Omega)$ is
any function which vanishes at infinity of $\Omega$, then clearly
the function $\Theta(x)V(x)$ is a (resp. semi-) small perturbation
of the operator $P$ (resp. $P^*$) in $\Omega$.

4. Suppose that $G^{\Omega}_{P}$ and $ G^{\Omega}_{P-|V|}$ are
equivalent (resp. semi-equivalent). Using the resolvent equation
it follows that the best equivalence (resp. semi-equivalence)
constants of $G^{\Omega}_{P}$ and $ G^{\Omega}_{P\pm |V^*_j|}$
tend to $1$ as $j \to \infty$ if and only if $V$ is a (resp.
semi-) small perturbation of $P$ (resp. $P^*$) in $\Omega$.
Therefore, zero-order perturbations of the type studied by
A.~Ancona in \cite{An97} provide us with a huge and almost optimal
class of examples of small perturbations. (see also \cite[and the
references therein]{Ai,M86,M97,P89}).
 }\end{Rems}
 A.~Grigor'yan and W.~Hansen \cite{GH}
have introduced the following notions of perturbations.
\begin{Def}\label{defhbig}{\em Let $P$ be a subcritical operator in
$\Omega$, and fix $h\in \mathcal{C}_P(\Omega)$. A nonnegative
function $V$ is called {\em $h$-big on $\Omega$} if any solution
$v$ of the equation $(P+V)v=0$ in $\Omega$ satisfying $0\leq v\leq
h$ is identically zero. $V$ is {\em non-$h$-big on $\Omega$} if
$V$ is not $h$-big on $\Omega$.
 }\end{Def}
\begin{remark}\label{remhbig}{\em If $V$ is $H$-integrable
perturbation of $P$, then it is non-$h$-big for any $h\in
\mathcal{C}_P(\Omega)$ (see Proposition~\ref{propinteg}).
 }\end{remark}
 The following notion of perturbation does not involve
Green functions.
\begin{Def} \label{weakpert}
{\rm Let $P$ be a subcritical operator in $\Omega\subseteq X$. A
function  $V \in C^{\ga}(\Omega)$ is said to be {\em a weak
perturbation} of the operator $P$ in $\Omega$ if the following
condition holds true.
\begin{description}
\item[($\ast$)] For every $\lambda \in \mathbb{R}$ there exists $N \in \Nat$ such
that the operator $P-\lambda V^*_n(x)$ is subcritical in $\Omega$
for any $n \geq N$.
\end{description}

A function $V \in C^{\ga}(\Omega)$ is said to be {\em a weak
perturbation} of a critical operator $P$ in $\Omega$ if there
exists a nonzero, nonnegative function $W \in C^{\ga}_0(\Omega)$
such that the function $V$ is a weak perturbation of the
subcritical operator $P+W$ in  $\Omega$. }\end{Def}
\begin{Rems}{\em
1. If $V$ is a weak perturbation of $P$ in $\Omega$, then
$\partial S=S_0$ and $\gl_\infty(P,\Omega,\pm V)=\infty$
(\cite{P98}, see also Theorem \ref{extthm}).

2. If $V$ is a semismall perturbation  of $P$ in $\Omega$, then
$|V|$ is a  weak perturbation of $P$ in $\Omega$, but $G$-bounded
perturbations are not necessarily weak.

3. Let $d\geq 3$. By the Cwikel-Lieb-Rozenblum estimate, if $V\in
L^{d/2}(\mathbb{R} ^d)$, then $|V|$ is a weak perturbation of
$-\Gd$ in $\mathbb{R} ^d$. On the other hand, $(1+|x|)^{-2}$ is
not a weak perturbation of $-\Gd$ in $\mathbb{R} ^d$, while for
any $\varepsilon
>0$ the function $(1+|x|)^{-(2+\varepsilon)}$ is a small perturbation
of $-\Gd$ in $\mathbb{R} ^d, d\geq 3$ \cite{M86,P88}. }\end{Rems}
\mysection{Indefinite weight}\label{secindef}
Consider the Schr\"{o}dinger operator $H_{\lambda}:= -\Delta
-\lambda W$ in $\mathbb{R} ^d$, where $\lambda \in \mathbb{R} $ is
a spectral parameter and $W\in C_0^{\infty}(\mathbb{R} ^d),
W\not\equiv 0$. Since $-\Delta$ is subcritical in $\mathbb{R} ^d$
if and only if $d \geq 3$, it follows that for $d \geq 3$ the
Schr\"{o}dinger operator $H_{\lambda}$ has no bound states
provided that $|\lambda|$ is sufficiently small. On the other
hand, for $d=1,2,$ B.~Simon proved the following sharp result.
\begin{Thm}[\cite{S76}] \label{Thmsimon}
Suppose that $d =1,2$, and let $W\in C_0^{\infty}(\mathbb{R} ^d),
W\not\equiv 0$. Then $H_{\lambda}=-\Delta -\lambda W$ has a
negative eigenvalue for all negative  $\lambda$ if and only if
$\int_{\mathbb{R} ^d} W(x) dx \leq 0.$
\end{Thm}
The following result extends Theorem~\ref{Thmsimon}  to the case
of a weak perturbation of a general critical operator in $\Omega$.
\begin{Thm}[\cite{P98}] \label{mainthmindef}
Let $P$ be a critical operator in $\Omega$, and $W \in
C^{\ga}(\Omega)$ a weak perturbation of the operator $P$ in
$\Omega$. Denote by $\varphi _{0}$ (resp. $\varphi^*_{0}$) the
ground state of the operator $P$ (resp. $P^*$) in $\Omega$ such
that $\varphi _{0}(x_0)=1$ (resp. $\varphi^*_{0}(x_0)=1$). Assume
that $W\varphi_0 \varphi^*_0\in L^1(\Omega)$.

(i) If there exists $\lambda<0$ such that $P-\lambda W(x)$ is
subcritical in $\Omega$, then \be \label{integralcon}
\int_{\Omega}W(x)\varphi_0(x) \varphi^*_0(x)\dx > 0. \ee

(ii) Assume that for some nonnegative, nonzero function $V \in
C^\ga _0(\Omega)$ there exists $\tilde{\lambda}<0$  and a positive
constant $C$ such that \be \label{cond} \Green{\Omega}{P+V-\lambda
W}{x}{x_0}\leq C\varphi_0(x) \quad \mbox{and} \quad
\Green{\Omega}{P+V-\lambda W}{x_0}{x}\leq C\varphi^*_0(x) \ee for
all $x\in \Omega \setminus \Omega _1$ and $\tilde{\lambda} \leq
\lambda <0$. If the integral condition (\ref{integralcon}) holds
true, then there exists $\lambda<0$ such that $P-\lambda W(x)$ is
subcritical in $\Omega$.

(iii) Suppose that $W$ is a semismall perturbation of the
operators $P+V$ and $P^*+V$ in $\Omega$, where $V\gvertneqq 0$, $V
\in C^\ga _0(\Omega)$ . Then there exists $\lambda<0$ such that
$P-\lambda W(x)$ is subcritical in $\Omega$ if and only if
\eqref{integralcon} holds true.
\end{Thm}
\mysection{Large time behavior of the heat kernel}\label{sechetk}
As was already mentioned in Section~\ref{secpreliminaries}, the
large time behavior of the heat kernel is closely related to
criticality (see for example Lemma~\ref{lemheatcrit}). In the
present section we elaborate this relation further more.

Suppose that $\lambda_0(P,\Omega,\mathbf{1})\geq 0$. We consider
the parabolic operator $L$
\begin{equation}\label{eqL}
  Lu=u_t+Pu \qquad \mbox{ on } \Omega\times (0,\infty).
\end{equation}  \noindent
 We denote by $\mathcal{H}_P(\Omega\times (a,b))$ the cone of all
nonnegative solutions of the equation $Lu=0$ in $\Omega\times
(a,b)$. Let $k_P^{\Omega}(x,y,t)$ be the heat kernel of the
parabolic operator $L$ in $\Omega$.

If $P$ is critical in $\Omega$, we denote by $\varphi_0$ the
ground state of $P$ in $\Omega$ satisfying $\varphi_0(x_0)=1$. The
corresponding ground state of $P^*$ is denoted by $\varphi^*_0$.

\begin{definition}{\em A critical operator $P$ is said to be {\em positive-critical}
in $\Omega$ if $\varphi_0\varphi_0^*\in L^1(\Omega)$, and {\em
null-critical} in $\Omega$ if $\varphi_0\varphi_0^*\not\in
L^1(\Omega)$. }
\end{definition}
\begin{theorem}[\cite{P92,P04}]\label{mainthmhk}
Suppose that $\lambda_0\geq 0$. Then for each  $x,y\in \Omega$
 \begin{equation*}\label{eqlimhk}
\lim_{t\to\infty} \mathrm{e}^{\lambda_0 t}k_P^{\Omega}(x,y,t)\!=\!
  \begin{cases}
    \dfrac{\varphi_0(x)\varphi_0^*(y)}{\int_\Omega\!
\varphi_0(z)\varphi_0^*(z)\,\mathrm{d}z} & \text{if }
P\!-\!\lambda_0 \text{ is positive-critical},
\\[5mm]
    0 & \text{otherwise}.
  \end{cases}
 \end{equation*}
Moreover, we have the following Abelian-Tauberian type relation
\begin{equation}\label{eqgreen}
\lim_{t\to\infty} \mathrm{e}^{\lambda_0 t}k_P^{\Omega}(x,y,t)=
\lim_{\lambda\nearrow\lambda_0}(\lambda_0-\lambda)\Green{\Omega}{P-\lambda}{x}{y}.
\end{equation}
\end{theorem}
\begin{remark}{\em
The first part of Theorem~\ref{mainthmhk} has been proved by
I.~Chavel and L.~Karp \cite{CK} in the {\em selfadjoint} case.
Later, B. Simon gave a shorter proof for the selfadjoint case
using the spectral theorem and elliptic regularity \cite{S93}.
 }\end{remark}

We next ask how fast $\lim_{t\to\infty} \mathrm{e}^{\lambda_0
t}k_P^{\Omega}(x,y,t)$ is approached.  It is natural to conjecture
that the limit is approached equally fast for different points
$x,y\in \Omega$. Note that in the context of Markov chains, such
an {\em (individual) strong ratio limit property} is in general
not true \cite{Chu}. The following conjecture was raised by
E.~B.~Davies \cite{Dheat} in the selfadjoint case.
\begin{conjecture}\label{conjD}
Let $Lu=u_t+P(x, \partial_x)u$ be a parabolic operator which is
defined on $\Omega\subseteq X$. Fix a reference point $x_0\in
\Omega$. Then
\begin{equation}\label{eqconjD}
\lim_{t\to\infty}\frac{k_P^\Omega(x,y,t)}{k_P^\Omega(x_0,x_0,t)}=a(x,y)
\end{equation}
exists and is positive for all $x,y\in \Omega$.
\end{conjecture}
If Conjecture~\ref{conjD} holds true, then for any fixed $y\in
\Omega$ the limit function $a(\cdot,y)$ is a positive solution of
the equation $(P-\lambda_0)u = 0$ which is (up to a multiplicative
function) a parabolic Martin function in
$\mathcal{H}_P(\Omega\times \mathbb{R}_-)$ associated with any
Martin sequence of the form $(y, t_n)$ where $t_n\to-\infty$
  (see \cite[and the references
therein]{Dheat,IWPT} for further partial results).
\mysection{Nonuniqueness of the positive Cauchy\\ problem and
intrinsic ultracontractivity }\label{secup}
In this section we discuss the uniqueness the Cauchy problem
\begin{equation}\label{eqL1}
  \begin{cases}
   Lu:=u_t+Pu=0  & \mbox{ on } \Omega\times (0,T), \\
    u(x,0)=u_0(x) & \mbox{ on } \Omega,
  \end{cases}
\end{equation}  \noindent
in the class of nonnegative continuous solutions. So, we always
assume that $u_0 \in C(X)$, and $u_0\geq 0$.
\begin{definition}\label{defup}{\em
A {\em solution of the positive Cauchy problem} in
$\Omega_T\!:=\!\Omega\times [0,T)$ with initial data $u_0$ is a
nonnegative continuous function in $\Omega_T$ satisfying
$u(x,0)=u_0(x)$, and $Lu=0$ in $\Omega\times(0,T)$ in the
classical sense.

We say that the {\em uniqueness of the positive Cauchy problem}
(UP) for the operator $L$ in $\Omega_T$ holds,  when any two
solutions of the positive Cauchy problem satisfying the same
initial condition are identically equal in $\Omega_T$.
 }\end{definition}

Let $u \in {\cal C}_{P}(\Omega)$. By the parabolic generalized
maximum principle, either \be \label{kueu} \int_\Omega\!\!
k(x,y,t)u(y)\!\dy\!=\! u(x)\;\;\mbox{for some (and hence for all)
$x \in \Omega,\, t>0$,} \ee or \be \label{kuslu} \int_\Omega\!\!
k(x,y,t)u(y)\!\dy\!<\!u(x)\;\;\mbox{for some (and hence for all)
$x \in \Omega,\, t>0,$} \ee see for example \cite{EBD}. Note that
both sides of \eqref{kuslu} are solutions of the positive Cauchy
problem \eqref{eqL1} with the same initial data $u_0=u$.
Therefore, in order to show that UP does not hold for the operator
$L$ in $\Omega$, it is sufficient to show that \eqref{kuslu} holds
true for some $u \in {\cal C}_{P}(\Omega)$. It is easy to show
\cite{EBD} that
 (\ref{kuslu}) holds true if and only if there exists
$\lambda <0$ such that \be \label{Gluslu}
-\lambda\int_\Omega\Green{\Omega}{P-\lambda}{x}{y}u(y)\dy < u(x)
\ee
for some (and hence for all) $x\in \Omega$. Furthermore, it
follows from \cite{M95} that (\ref{Gluslu}) is satisfied if
\be \label{Gufnt} \int_\Omega\Green{\Omega}{P}{x}{y}u(y)\dy <
\infty \ee for some (and hence for all) $x\in \Omega$. Thus, we
have:
\begin{corollary}\label{cornup}
If $\mathbf{1}$ is an $H$-integrable perturbation of a subcritical
operator $P$  in $\Omega$, then the positive Cauchy problem is not
uniquely solvable.
\end{corollary}
\begin{remarks}{\em
1. A positive solution $u \in {\cal C}_{P}(\Omega)$ which
satisfies (\ref{kueu}) is called a {\em positive invariant
solution}.  If $P\mathbf{1}=0$ and (\ref{kueu}) holds for
$u=\mathbf{1}$ one says that {\em $L$ conserves probability in}
$\Omega$ (see \cite{G99}). We note that if $P$ is critical, then
the ground state $\varphi_0$ is a positive invariant solution. It
turns out that there exists a complete Riemannian manifold $X$
which does not admit any positive invariant harmonic function,
while $\lambda_0(-\Delta,X,\mathbf{1})=0$ \cite{PStroock}.

2. For necessary and sufficient conditions for UP, see
\cite{IM,M05} and the references therein.
 }\end{remarks}


The following important notion was introduced by E.~B.~Davies and
B.~Simon for Schr\"odinger operators \cite{DS84,DS86a,DS86b}.
\begin{definition}{\em
Suppose that $P$ is symmetric. The Schr\"odinger semigroup
$\mathrm{e}^{-tP}$ associated with the heat kernel
$k_P^{\Omega}(x,y,t)$ is called {\em intrinsic ultracontractive}
(IU) if $P-\lambda_0$ is positive-critical in $\Omega$ with a
ground state $\varphi_0$, and for each $t>0$ there exists a
positive constant $C_t$ such that
$$C_t^{-1}\varphi_0(x)\varphi_0(y)\leq k_P^{\Omega}(x,y,t)\leq C_t\varphi_0(x)\varphi_0(y)
\qquad \forall x,y\in \Omega.$$ }\end{definition} .

\vspace{-8mm}

\begin{remarks}\label{remIU} {\em
1. If $\mathrm{e}^{-tP}$  is IU, then
\begin{equation}\label{eqlimhk1}
\lim_{t\to\infty} \mathrm{e}^{\lambda_0 t}k_P^{\Omega}(x,y,t)=
    \dfrac{\varphi_0(x)\varphi_0(y)}{\int_\Omega
[\varphi_0(z)]^2\,\mathrm{d}z}
 \end{equation}
{\em uniformly} in $\Omega\times \Omega$ (see for example
\cite{Ba}, cf. Theorem~\ref{mainthmhk}).

2. If $\Omega$ is a bounded uniformly H\"older domain of order
$0<\alpha<2$, then $\mathrm{e}^{-t(-\Delta)}$  is IU on $\Omega$
\cite{Ba}.

3. Let $\alpha\geq 0$. Then $\mathrm{e}^{-t(-\Delta+|x|^\alpha)}$
is IU on $\mathbb{R}^d$ if and only if $\alpha> 2$.
  }\end{remarks}

Intrinsic ultracontractivity is closely related to perturbation
theory of positive solutions and hence to UP, as the following
recent result of M.~Murata and M.~Tomisaki demonstrates.
\begin{theorem}[\cite{M97,MT2006}]\label{thmiusp}
Suppose that $P$ is a subcritical symmetric operator, and that the
Schr\"odinger semigroup $\mathrm{e}^{-tP}$ is IU on $\Omega$. Then
${\bf 1}$ is a small perturbation of $P$ on $\Omega$. In
particular, UP does not hold in $\Omega$.
\end{theorem}
On the other hand, there are  planner domains such that ${\bf 1}$
is a small perturbation of the Laplacian, but the semigroup
$\mathrm{e}^{-t(-\Delta)}$ is not IU (see \cite{BD} and
\cite{P99}).
\mysection{Asymptotic behavior of eigenfunctions}\label{seceigen}
In this section, we assume that $P$ is symmetric and discuss
relationships between perturbation theory, Martin boundary, and
the asymptotic behavior of weighted eigenfunctions in some general
cases (for other relationships between positivity and decay of
Schr\"odinger eigenfunctions see, \cite{Agmon84,S82,S00}).
\begin{Thm} \label{extthm}
(i) Let $V\in C^\ga(\Omega)$ be a positive function. Suppose that
$P$ is a symmetric, nonnegative operator on $L^2(\Omega,V(x)dx)$
with a domain $C_0^\infty(\Omega)$. Assume that $V$ is a weak
perturbation of the operator $P$ in $\Omega$. suppose that $P$
admits a (Dirichlet) selfadjoint realization $\tilde{P}$ on
$L^2(\Omega,V(x)dx)$. Then $\tilde{P}$ has a purely discrete
nonnegative spectrum (that is,
$\gs_{\mathrm{ess}}(\tilde{P})=\emptyset$). Moreover,
$$\gs(\tilde{P})=\gs_{\mathrm{discrete}}(\tilde{P})=\gs_{\mathrm{point}}(\tilde{P})=
\{\gl_n\}_{n=0}^\infty,$$ where $\lim_{n\to \infty}\gl_n=\infty$.
In particular, if $\gl_0:=\gl_0(P,\Omega,V)>0$, then the natural
embedding $E:{\cal H} \longrightarrow L^2(\Omega,V(x)dx)$ is
compact, where ${\cal H}$ is the completion of
$C_0^\infty(\Omega)$  with respect to the inner product induced by
the corresponding quadratic form.

(ii) Assume further that $P$ is subcritical and $V$ is a semismall
perturbation of the  operator $P$ in $\Omega$. Let
$\{\varphi_n\}_{n=0}^\infty$ be the set of the corresponding
eigenfunctions ($P\varphi_n=\gl_nV\varphi_n$). Then for every
$n\geq 1$ there exists a positive constant $C_n$ such that
\be\label{efest} |\varphi_n(x)|\leq C_n \varphi_0(x). \ee

(iii) For every $n\geq 1$,  the function $\varphi_n/\varphi_0$ has
a continuous extension $\psi_n$ up to the Martin boundary
$\partial_{P}^M\Omega$, and $\psi_n$ satisfies
$$\psi_n(\xi)\!=\!
(\psi_0(\xi))^{-1}\!\gl_n\!\!\!\int_{\Omega}\!\!K^\Omega_{P}(z,\xi)V(z)\varphi_n(z)\!\dz\!=\!
\frac{\gl_n\!\int_{\Omega}\!K^\Omega_{P}(z,\xi)V(z)\varphi_n(z)\!\dz}
{\gl_0\!\int_{\Omega}\!K^\Omega_{P}(z,\xi)V(z)\varphi_0(z)\!\dz}
$$
for every $\xi\in\partial_{P}^M\Omega$, where $\psi_0$ is the
continuous extension of $\varphi_0/\Green{\Omega}{P}{\cdot}{x_0}$
to the Martin boundary $\partial_{P}^M\Omega$.
\end{Thm}
\begin{Rems}\label{IUrem}{\em
1. By \cite{DS84}, the semigroup $\mathrm{e}^{-t\tilde{P}}$ is IU
if and only if the pointwise eigenfunction estimate (\ref{efest})
holds true with $C_n=c_t\exp (t\gl_n)\|\varphi_n\|_2$, for every
$t>0$ and $n\geq 1$. Here $c_t$ is a positive function of $t$
which may be taken as the function such that
$k_P^\Omega(x,y,t)\leq c_t\varphi_0(x)\varphi_0(y)$, where
$k_P^\Omega$ is the corresponding heat kernel. It follows that if
$\mathrm{e}^{-t\tilde{P}}$ is IU, then the pointwise eigenfunction
estimate (\ref{efest}) holds true with $C_n=\inf_{t>0}\{c_t\exp
(t\gl_n)\}\|\varphi_n\|_2$. We note that in general $\{C_n\}$ is
unbounded \cite{GD}.

Recall that if $\mathrm{e}^{-t\tilde{P}}$ is IU, then $\mathbf{1}$
is a small perturbation of $P$ (see Theorem~\ref{thmiusp}). In
particular, part (iii) of Theorem~\ref{extthm} implies that if
$\mathrm{e}^{-t\tilde{P}}$ is IU, then for any $n\geq 1$, the
quotient $\varphi_n/\varphi_0$ has a continuous extension $\psi_n$
up to the Martin boundary $\partial_{P}^M\Omega$.

2. M.~Murata \cite{Mu90} proved part (ii) of Theorem \ref{extthm}
for the special case of bounded Lipschitz domains. See also
\cite{HO} for related results on the asymptotic behavior of
eigenfunctions of Schr\"{o}dinger operators in $\mathbb{R}^d$.

 }\end{Rems}
\mysection{Localization of binding}\label{seclocalization}
Let $V\in C^\alpha(\mathbb{R}^d)$ and $R\in \mathbb{R}^d$,
throughout this section we use the notation $V^R(x):=V(x-R)$. For
$j=1,2$, let $V_{j}$ be small perturbations of the Laplacian in
$\mathbb{R}^d, d\geq 3$, and assume that the operators
$P_{j}:=-\Delta + V_{j}(x)$ are nonnegative on $\core$. We
consider the Schr\"{o}dinger operator \be \label{dePR} P_{R} :=
-\Delta + V_{1}(x) + V_{2}^R(x) \ee defined on $\mathbb{R}^d$, and
its ground state energy
$E(R):=\lambda_0(P_R,\mathbb{R}^d,\mathbf{1})$. In this section we
discuss the asymptotic behavior of $E(R)$ as $\abs{R}\to\infty$, a
problem which was studied by M.~Klaus and B.~Simon in
\cite{KS79,S80} (see also \cite{Plocal,Pins95}). The motivation
for studying the asymptotic behavior of  $E(R)$ comes from a
remarkable phenomenon known as the Efimov effect for a three-body
Schr\"{o}dinger operator (for more details, see for example
\cite{T1}).
\begin{definition} \label{Kato}
Let $d\geq 3$. The space of functions \be \label{Kndef}
K_{d}^{\infty}:=\left\{ V \in C^{\alpha}(\mathbb{R}^d) |
\lim_{M\rightarrow
\infty}\sup_{x\in\mathbb{R}^d}\int_{|z|>M}\frac{|V(z)|}
{|x-z|^{d-2}}\dz=0 \right\} \ee is  called  the  Kato  class  at
infinity.
\end{definition}
\begin{remark}\label{remkato}
{\em Let $d\geq 3$. If $V\in K_{d}^{\infty}$, then  $V$ is a small
perturbation of the Laplacian in $\mathbb{R}^d$.
 }\end{remark}

\begin{theorem}[\cite{Plocal}] \label{Thm1}
Let $d\geq 3$. For $j=1,2$, let $V_{j}(x)\in K_{d}^{\infty}$ be
two functions such that the operators $P_{j}=-\Delta + V_{j}(x)$
are subcritical in $\mathbb{R}^d$. Then there exists $r_{0}>0$
such that the operator $P_{R}$ is subcritical for any $R\in
\mathbb{R}^d \setminus B(0,r_0)$. In particular, $E(R)=0$ for all
$\abs{R}\geq r_0$.
\end{theorem}
Assume  now that  the  operators  $P_{j}=-\Delta + V_{j}(x)
,\;j=1,2$, are  {\em critical} in  $\mathbb{R}^d$. It turns out
that in this case, there exists $r_{0}>0$ such that $E(R)<0$ for
$\abs{R}\geq r_0$, but the asymptotic behavior of $E(R)$ depends
on the dimension $d$, as the following theorems demonstrate (cf.
\cite[the remarks in pp. 84 and 87]{KS79}).
\begin{theorem}[\cite{T1}] \label{Thmlob3}
Let $d=3$. Assume that the potentials $V_{j}, \;j=1,2$ satisfy $
|V_{j}(x)| \leq C\langle x \rangle^{-\beta}$ on  $\mathbb{R}^3$,
where $\langle x \rangle:=(1+|x|^2)^{1/2}$, $\beta
>2$, and $C>0$. Suppose that
$P_j=-\Delta + V_{j}(x)$ is critical in $\mathbb{R}^3$ for
$j=1,2$.

Then there exists $r_{0}>0$ such that the operator $P_{R}$ is
supercritical for any $R\in \mathbb{R}^3 \setminus B(0,r_0)$.
Moreover, $E(R)$ satisfies
  \be \label{R2ER} \lim_{\abs{R} \rightarrow
\infty}\abs{R}^{2}E(R)=-\beta^2 <-1/4, \ee where $\beta$ is the
unique root of the equation $s=\mathrm{e}^{-s}$.
\end{theorem}

\begin{theorem}[\cite{P96}] \label{Thmlob4}
Let $d=4$. Assume that for $j=1,2$ the operators $P_{j}=-\Delta +
V_{j}(x)$ are critical in $\mathbb{R}^4$, where $V_{j}\in
C^\alpha_0(\mathbb{R}^4)$.

Then there exists $r_{0}>0$ such that the operator $P_{R}$ is
supercritical for any  $R\in \mathbb{R}^4 \setminus B(0,r_0)$.
Moreover, there exists a  positive constant $C$ such that $E(R)$
satisfies \be \label{ER6est} -C\abs{R}^{-2} \leq E(R) \leq
-C^{-1}|R|^{-2}(\log |R|)^{-1} \;\; \ \mbox{for all $|R|\geq
r_{0}$.} \ee
\end{theorem}
\begin{theorem}[\cite{Plocal}] \label{Thmlob5}
Let $d\geq 5$. Suppose that $V_{j}, \;j=1,2$ satisfy $ |V_{j}(x)|
\leq C\langle x \rangle^{-\beta}$ in $\mathbb{R}^d$, where $\beta
> d-2$, and $C>0$. Assume that the operators
$P_{j}=-\Delta + V_{j}(x),\;j=1,2$, are critical in
$\mathbb{R}^d$.

Then there exists $r_{0}>0$ such that the operator $P_{R}$ is
supercritical for any $R\in \mathbb{R}^d \setminus B(0,r_0)$.
Moreover, there exists a  positive constant $C$ such that $E(R)$
satisfies \be \label{ER5est} -C|R|^{2-d} \leq E(R) \leq
-C^{-1}|R|^{2-d} \;\; \ \mbox{for all $|R|\geq r_{0}$.} \ee
\end{theorem}
What distinguishes $d\geq 5$ from $d=3,4$, is that for a
short-range potential $V$, the ground state of a critical operator
$-\Delta + V(x)$ in $\mathbb{R}^d$ is in $L^{2}(\mathbb{R}^d)$ if
and only if $d\geq 5$ (see \cite{S81} and Theorem~\ref{thmssp}).
\mysection{The shuttle operator}\label{secshttle}
In this section we present an intrinsic criterion which
distinguishes between subcriticality, criticality and
supercriticality of the operator $P$ in $\Omega$. This criterion
depends only on the norm of a certain linear operator $S$, called
the {\em shuttle operator} which is defined on $C(\partial D)$,
where $D\Subset\Omega$.

The shuttle operator was introduced for Schr\"odinger operators on
$\mathbb{R}^d$ in \cite{C,CV,Z1,Z2}. Using Feynman-Kac-type
formulas \cite{Spath}, F.~Gesztesy and Z.~Zhao \cite{GZ2,Z2} have
studied the shuttle operator for Schr\"{o}dinger operators in
$\mathbb{R}^d$ with short-range potentials (see also \cite{GZ1}),
and its relation to the following problem posed by B.~Simon.
\begin{problem}[\cite{S81,S82}]\label{problem}
Let $V\in L^2_{\mathrm{loc}}(\mathbb{R}^2)$. Show that if the
equation $(-\Delta\!+\!V)u\!=\!0$ on $\mathbb{R}^2$ admits a
positive $L^\infty$-solution, then $-\Delta\!+\!V$ is critical.
\end{problem}
Gesztesy and Zhao used the shuttle operator and proved that for
{\em short-range} potentials on $\mathbb{R}^2$, the above
condition is a necessary and sufficient condition for criticality
(see also \cite{M84} and Theorem~\ref{thmssp} for similar results,
and Theorem~\ref{thmDKS} for the complete solution). On the other
hand, Gesztesy and Zhao showed in \cite[Example~4.6]{GZ1} that
there is a critical Schr\"odinger operator on $\mathbb{R}$ with
`almost' short-range potential such that its ground state behaves
logarithmically.

Let $P$ be an elliptic operator of the form (\ref{P}) which is
defined on $\Omega$. We assume that the following assumption {\bf
(A)} holds:
\begin{description}
\item [(A)]
There exist four smooth, relatively compact subdomains
$\Omega_j,\; 0\le j\le 3$, such that $\overline{\Omega_{j}}\subset
\Omega_{j+1},\;j=0,1,2$, and such that
$\mathcal{C}_{P}(\Omega_3)\neq \emptyset$ and
$\mathcal{C}_{P}(\Omega_{0}^*)\neq \emptyset$.
\end{description}
\begin{remarks} \label{rem(A)}
{\em 1. If assumption {\bf (A)} is not satisfied, then we shall
say that the spectral radius of the shuttle operator is infinity.
In this case, it is clear that $P$ is supercritical in $\Omega$.

2. Assumption {\bf (A)} does not imply that
$\mathcal{C}_{P}(\Omega)\neq \emptyset$. }
\end{remarks}

Fix an exhaustion $\{\Omega_{j}\}_{j=0}^{\infty}$ of $\Omega$,
such that $\Omega_j$  satisfy assumption {\bf (A)} for $0\le j\le
3$. By assumption {\bf (A)} the Dirichlet problem \bea Pu=0
\;\;\;\mbox{in $\Omega_2$},\quad u=f  \;\;\;&\mbox{on $\partial
\Omega_2$} \eea is uniquely solved in $\Omega_2$ for any $f \in
C(\partial \Omega_2)$, and we denote the corresponding operator
from $C(\partial \Omega_2)$ into $C(\Omega_2)$ by $T_{\Omega_2}$.
Moreover, for every $f \in C(\partial \Omega_1)$, one can uniquely
solve the exterior Dirichlet problem in the outer domain
$\Omega_1^*$, with `zero' boundary condition at infinity of
$\Omega$. So, we have an operator $T_{\Omega_1^*}: C(\partial
\Omega_1) \to C(\Omega_1^*)$ defined by
$$T_{\Omega_1^*}f(x):=\lim_{j\to\infty}u_{f,j}(x),$$
where $u_{f,j}$ is the solution of the Dirichlet boundary value
problem:
$$
Pu=0 \;\;\;\mbox{in $\Omega_1^*\cap \Omega_j$},\quad u=f
\;\;\;\mbox{on $\partial \Omega_1^*$},\quad u=0  \;\;\;\mbox{on
$\partial (\Omega_1^*\cap \Omega_j) \setminus
\partial \Omega_1^*$}.
$$
For any  open set $D$ and $F \Subset D$, we denote by $R^D_{F}$
the restriction map $f\mapsto f\!\!\mid_F$ from $C(D)$ into
$C(F)$. The {\em shuttle operator} $S:C(\partial \Omega_1)
\longrightarrow C(\partial \Omega_1)$ is defined as follows: \be
\label{defshuttle} S:=R^{\Omega_2}_{\partial
\Omega_1}T_{\Omega_2}R^{\Omega_1^*}_{\partial
\Omega_2}T_{\Omega_1^*}\,. \ee We denote the spectral radius of
the operator $S$ by $r(S)$. We have
\begin{Thm}[\cite{Pshuttle}] \label{mainthmshuttle}  The operator $P$ is
subcritical, critical, or supercritical in $\Omega$ according to
whether $r(S)<1$, $r(S)=1$, or  $r(S)>1$.
\end{Thm}
The proof of Theorem~\ref{mainthmshuttle} in \cite{Pshuttle} is
purely analytic and relies on the observation that (in the
nontrivial case) $S$ is  a positive compact operator defined on
the Banach space $C(\partial \Omega_1)$. Therefore, the
Krein-Rutman theorem implies that there exists a simple principal
eigenvalue $\nu_0>0$, which is equal to the norm (and also to the
spectral radius) of $S$, and that the corresponding principal
eigenfunction is strictly positive. It turns out, that the
generalized maximum principle holds in any smooth subdomain
$D\Subset \Omega$ if and only if $\nu_0 \le 1$, and that $\nu_0
<1$ if and only if $P$ admits a positive minimal Green function in
$\Omega$.

The shuttle operator can be used to prove localization of binding
for certain {\em nonselfadjoint} critical operators (see
\cite{Pshuttle}).
\mysection{Periodic operators}\label{secperiod}
In this section we restrict the form of the operator. Namely, we
assume that $P$ is defined on $\mathbb{R}^d$ and that the
coefficients of $P$ are $\mathbb{Z}^d$-periodic. For such
operators, we introduce a function $\Lambda$ that plays a crucial
role in our considerations. Its properties were studied in detail
in \cite{A1a,KS87,LP,MT,Pinsper}. Consider the function $\Lambda
:\mathbb{R}^d\rightarrow \mathbb{R}$ defined by the condition that
the equation $Pu=\Lambda (\xi )u$ on $\mathbb{R}^d$ has a positive
{\em Bloch solution} of the form
\begin{equation}\label{positiveBloch}
u_{\,\xi }(x)=e^{\xi \cdot x}\varphi_{\,\xi }(x),
\end{equation}
where $\xi\in\mathbb{R}^d$, and $\varphi_{\,\xi }$ is a positive
$\mathbb{Z}^d$-periodic function.

\begin{theorem}
\label{Lambda-lemma}
\begin{enumerate}
\item  The value $\Lambda (\xi )$ is uniquely determined for any $\xi \!\in\!
\mathbb{R}^d$.

\item  The function $\Lambda$ is bounded from above, strictly
concave, analytic, and has a nonzero gradient for any $\xi \in
\mathbb{R}^d$  except at its maximum point.

\item  For $\xi \in
\mathbb{R}^d$, consider the operator $P(\xi ):=e^{-\xi\cdot
x}Pe^{\xi\cdot x}$ on the torus $\mathbb{T}^d$.  Then $\Lambda
(\xi )$ is the principal eigenvalue of $P(\xi )$ with  a positive
eigenfunction $\varphi_{\,\xi }$. Moreover, $\Lambda (\xi )$ is
algebraically simple.

\item The Hessian of $\Lambda (\xi )$ is nondegenerate at all points $\xi \in
\mathbb{R}^d$.
%
\end{enumerate}
\end{theorem}
Let us denote
\begin{equation}
\Lambda_{0} =\max_{\xi \in \mathbb{R}^d}\Lambda (\xi ).
\label{Lambda}
\end{equation}
It follows from \cite{A1a,LP,Pinsper} that
$\Lambda_{0}=\lambda_0$, and that $P-\Lambda_{0}$ is critical if
and only if $d=1,2$ (see also Corollary~\ref{corper}). Thus, in
the self-adjoint case, $\Lambda_{0}$ coincides with the bottom of
the spectrum of the operator $P$.  Assume that $\Lambda_{0}\geq
0$. Then Theorem~\ref{Lambda-lemma} implies that the zero level
set
\begin{equation}
\Xi =\left\{ \xi \in \mathbb{R}^d|\;\Lambda (\xi )=0\right\}
\label{ksi}
\end{equation}
is either a strictly convex compact analytic surface in
$\mathbb{R}^d$ of dimension $d-1$ (this is the case if and only if
$\Lambda_{0}> 0$), or a singleton  (this is the case if and only
if $\Lambda_{0}=0$).

In a recent paper \cite{MT}, M.~Murata and T.~Tsuchida have
studied the exact asymptotic behavior at infinity of the positive
minimal Green function and the Martin boundary of such periodic
elliptic operators on $\mathbb{R}^d$.

Suppose that $\Lambda_0=\Lambda(\xi_0)>0$. Then $P$ is
subcritical, and for each $s$ in the unit sphere $S^{d-1}$ there
exists a unique $\xi_s\in \Xi$ such that
$$\xi_s\cdot s=\sup_{\xi\in\Xi}\,\{\xi\cdot s\}.$$ For $s\!\in\!
S^{d-1}$ take an orthonormal basis of $\mathbb{R}^d$ of the form
$\{e_{s,1},\ldots,e_{s,d-1},s\}$. For $\xi\in \mathbb{R}^d$, let
$\varphi_\xi$ and $\varphi^*_\xi$ be periodic positive solutions
of the equation $P(\xi)u=\Lambda(\xi)u$ and
$P^*(\xi)u=\Lambda(\xi)u$ on $\mathbb{T}^d$, respectively, such
that
$$\int_{\mathbb{T}^d}\varphi_{\xi}(x)\varphi^*_{\xi}(x)\mathrm\,{d}x=1.$$

\begin{theorem}[\cite{MT}]\label{thmmuratatsuchida1} 1. Suppose that $\Lambda_0>0$.
Then the minimal Green
function $G^{\mathbb{R}^d}_P$ of $P$ on $\mathbb{R}^d$ has the
following asymptotics as $\abs{x-y}\to\infty$:
$$G^{\mathbb{R}^d}_P(x,y)\!=\!\frac{\abs{\nabla\Lambda(\xi_s)}^{(d-3)/2}
\,\mathrm{e}^{-(x-y)\cdot\xi_s}\varphi_{\xi_s}(x)\varphi^*_{\xi_s}(y)}
{(2\pi\abs{x\!-\!y})^{(d-1)/2}[\mathrm{det}(-e_{s,j}\!\cdot\!
\Hess
\Lambda(\xi_s)e_{s,k})]^{1/2}}\left[1\!+\!O(\abs{x\!-\!y}^{-1})\right]\!,$$
where $s:=(x-y)/\abs{x-y}$.

2. Suppose that $\Lambda_0=\Lambda(\xi_0)=0$ and $d\geq 3$. Then
the minimal Green function $G^{\mathbb{R}^d}_P$ of $P$ on
$\mathbb{R}^d$ has the following asymptotics as
$\abs{x-y}\to\infty$:
$$G^{\mathbb{R}^d}_P(x,y)\!\!=\!\!\!\frac{}{}
\frac{2^{-1}\pi^{-d/2}\Gamma(\frac{d-2}{2})\,\mathrm{e}^{-(x-y)\cdot\xi_0}\varphi_{\xi_0}(x)\varphi^*_{\xi_0}(y)}
{\{\mathrm{det}[\Hess \Lambda(\xi_0)]\}^{1/2}\!\abs{[-\Hess
\Lambda(\xi_0)]^{-1/2}(x\!-\!y)}^{d-2}}
\!\!\left[1\!+\!O(\abs{x\!-\!y}^{-1})\right]\!\!.$$
\end{theorem}

Combining the results in \cite{A1a,MT}, we have the following
Martin representation theorem.
\begin{theorem}[\cite{A1a,MT}]\label{thmperrep}
Let $\Xi$ be the set of all $\xi \in \mathbb{R}^d$ such that the
equation $Pu=0$ admits a positive Bloch solution
$u_{\,\xi}(x)=e^{\xi\cdot x}\varphi_{\,\xi }(x)$ with
$\varphi_{\,\xi }(0)=1$. Then $u$ is a positive Bloch solution if
and only if $u$ is a minimal Martin function of the equation
$Pu=0$ in $\mathbb{R}^d$. Moreover, all Martin functions are
minimal. Furthermore, $u\in \mathcal{C}_P(\mathbb{R}^d)$ if and
only if there exists a positive finite measure $\gm$ on $\Xi$ such
that
\[
u(x)=\int_\Xi u_{\,\xi }(x)\,\mathrm{d}\gm(\xi ).
\]
\end{theorem}
Theorem~\ref{thmperrep} (except the result that all Martin
functions are minimal) was extended by V.~Lin and the author to a
manifold with a group action \cite{LP}. It is assumed that $X$ is
a noncompact manifold equipped with an action of a group $G$ such
that $GV=X$ for a compact subset $V\Subset X$, and that the
operator $P$ is a $G$-invariant operator on $X$ of the form
\eqref{P}. If $G$ is finitely generated, then the set of all
normalized positive solutions of the equation $Pu=0$ in $X$ which
are also eigenfunctions of the $G$-action is a real analytic
submanifold $\Xi$ in an appropriate finite-dimensional vector
space $\mathcal{H}$. Moreover, if $\Xi$  is not a singleton, then
it is the boundary of a strictly convex body in $\mathcal{H}$. If
the group $G$ is {\em nilpotent}, then any positive solution in
${\cal C}_{P}(X)$ can be uniquely represented as an integral of
solutions over $\Xi$. In particular, $u\in {\cal C}_{P}(X)$ is a
positive minimal solution if and only if it is a positive solution
which is also an eigenfunction of the $G$-action.
\mysection{Liouville theorems for Schr\"odinger operators and
Criticality}\label{seccritliouville}
The existence and nonexistence of nontrivial bounded solutions of
the equation $Pu=0$ are closely related to criticality theory as
the following results demonstrate (see also
Section~\ref{secliouville}).
\begin{Pro}[\cite{G},{\cite[Lemma 3.4]{P99}}]\label{propinteg}
Suppose that $V$ is a nonzero, nonnegative function such that $V$
is an $H$-integrable perturbation of a subcritical operator $P$ in
$\Gw$ and let $u\in {\cal C}_{P}(\Gw)$. Then for any
$\varepsilon>0$ there exists $u_\varepsilon \in {\cal
C}_{P+\varepsilon V}(\Gw)$ which satisfies $0<u_\varepsilon\leq u$
and the resolvent equation \be \label{usubteq}
u_\varepsilon(x)=u(x) -
\varepsilon\int_{\Omega}\Green{\Omega}{P+\varepsilon
V}{x}{z}V(z)u(z)\,\dz. \ee In particular, if $P\mathbf{1}=0$, then
for any $\varepsilon>0$ the operator $P+\varepsilon V$ admits a
nonzero bounded solution.
\end{Pro}
 In \cite[Theorem~5]{DKS}, D.~Damanik, R.~Killip, and B.~Simon
proved a result which, formulated in the following new way,
reveals a complete answer to Problem~\ref{problem} posed by
B.~Simon in \cite{S81,S82} (see also \cite{GZ2,M84} and
Theorem~\ref{thmssp}). An alternative proof based on criticality
theory is presented below.
\begin{theorem}[\cite{DKS}]\label{thmDKS}
Let $d = 1$ or $2$, and $q \in L^2_{\mathrm{loc}}(\mathbb{R}^d)$.
Suppose that $H_q:=-\Delta+q$ has a bounded positive solution in
$\mathcal{C}_{H_q}(\mathbb{R}^d)$. If $V\in
L^2_{\mathrm{loc}}(\mathbb{R}^d)$ and both $H_{q\pm V}\geq 0$,
then $V=0$. In other words, $H_q$ is critical.
\end{theorem}
\begin{proof}
Theorem~\ref{thmconv} implies that we should indeed show that
$H_q$ is critical. Assume that $H_q$ is subcritical. Take a
nonzero nonnegative $W$ with a compact support. Then by
Theorem~\ref{thmssp}, there exists $\varepsilon>0$ such that
$H_{q-\varepsilon W}\geq 0$. Let $M<N$. For  $d=1$ take the cutoff
function
 \bean a_{M,N}(x):=
  \begin{cases}
  0 & |x|>N, \\
1 & |x|\leq M,\\
  1-\frac{|x|-M}{N-M} & M< |x|\leq N,
  \end{cases}
 \eean
and for $d=2$
 \bean a_{M,N}(x):=
\begin{cases}
  0 & |x|>N, \\
1& |x|\leq M,\\
  \frac{\log N-\log |x|}{\log N-\log M} & M< |x|\leq N.
  \end{cases}
 \eean
Let $\psi$ be a positive bounded solution of the equation $H_qu=0$
in $\mathbb{R}^d$. Then for appropriate $N, M$ with $M,N\to\infty$
(see \cite{DKS}), we have
\begin{multline*}\label{E}
  0<c<\varepsilon \int_{\mathbb{R}^d} W (a_{M,N} \psi)^2\dx \leq
 \int_{\mathbb{R}^d} \left[|\nabla (a_{M,N}\psi)|^2 + q (a_{M,N} \psi)^2\right]
  \dx =\\
  \int_{\mathbb{R}^d} |\nabla a_{M,N}|^2 \psi^2\dx\to
  0,
\end{multline*}
and this is a contradiction.
\end{proof}
\begin{remarks}\label{remlioiva}{\em
1. Theorem~\ref{thmDKS} is related to Theorem~1.7 in \cite{BCN}
which claims that for $d=1,2$, if $H_q$ admits a bounded solution
that changes its sign, then $\lambda_0\!<\!0$. This claim and
Theorem~\ref{thmDKS} do not hold for $d\!\geq\! 3$ \cite{B}.

2. For other relationships between perturbation theory of positive
solutions and Liouville theorem see \cite{G99,GH}.
 }\end{remarks}
After submitting the first version of the present article to the
editors, we proved the following result which generalized
Theorem~\ref{thmDKS} and the Liouville type theorems in
\cite{BCN}.
\begin{theorem}[\cite{P06}]\label{mainthmliouville} Let $\Omega\subset X$ be a domain.
Consider two Schr\"odinger operators defined on $\Omega$ of the
form
\begin{equation}\label{eqpj}
P_j:=-\nabla\cdot(A_j\nabla)+V_j\qquad j=0,1,
\end{equation}
such that  $V_j\in L^{p}_{\mathrm{loc}}(\Omega;\mathbb{R})$ for
some $p>{d}/{2}$, and $A_j:\Omega \rightarrow \mathbb{R}^{d^2}$
are measurable matrix valued functions such that for any $K\Subset
\Omega$  there exists $\mu_K>1$ such that \be \label{stell}
\mu_K^{-1}I_d\le A_j(x)\le \mu_K I_d \qquad \forall x\in K,
 \ee
where $I_d$ is the $d$-dimensional identity matrix.

 Assume that the following assumptions hold true.

\begin{itemize}
\item[(i)] The operator  $P_1$ is critical in $\Omega$. Denote
by $\varphi\in \mathcal{C}_{P_1}(\Omega)$ its ground state.

\item[(ii)]  $\lambda_0(P_0,\Omega,\mathbf{1})\geq 0$, and there exists a
real function $\psi\in H^1_{\mathrm{loc}}(\Omega)$ such that
$\psi_+\neq 0$, and $P_0\psi \leq 0$ in $\Omega$, where
$u_+(x):=\max\{0, u(x)\}$.

\item[(iii)] The following matrix inequality holds
\begin{equation}\label{psialephia}
\psi^2(x) A_0(x)\leq C\varphi^2(x) A_1(x)\qquad  \mbox{ a. e. in }
\Omega,
\end{equation}
where $C>0$ is a positive constant.
\end{itemize}
Then the operator $P_0$ is critical in $\Omega$, and $\psi$ is its
ground state. In particular, $\dim \mathcal{C}_{P_0}(\Omega)=1$
and $\lambda_0(P_0,\Omega,\mathbf{1})=0$.
\end{theorem}
The proof of Theorem~\ref{mainthmliouville} relies on
Theorem~\ref{mainky3}.
\begin{corollary}[\cite{Pinsper}]\label{corper}
Assume that the coefficients of the elliptic operator
$P:=-\nabla\cdot(A\nabla)+V$ are $\mathbb{Z}^d$-periodic on
$\mathbb{R}^d$. Then the operator $P-\lambda_0$ is critical in
$\mathbb{R}^d$ if and only if $d\leq 2$.
\end{corollary}
\begin{remark}\label{remcorper}{\em
One can use \cite{LP} to extend Corollary~\ref{corper}  to the
case of equivariant Schr\"odinger operators on cocompact
coverings. Let $X$ be a noncompact nilpotent covering of a compact
Riemannian manifold. Suppose that $P:=-\Delta+V$ is an equivariant
operator on $X$ with respect to its {\em nilpotent} deck group
$G$. Then $P-\lambda_0$ is critical in $X$ if and only if $G$ has
a normal subgroup of finite index isomorphic to $\mathbb{Z}^d$ for
$d\leq 2$.
 }\end{remark}
 \mysection{Polynomially growing solutions and Liouville Theorems}\label{secliouville}
 Let
$H=-\Delta +V$ be a Schr\"odinger operator on $\mathbb{R}^d$. Then
\v{S}nol's theorem asserts that, under some assumptions on the
potential $V$, if $H$ admits a polynomially growing solution of
the equation $Hu=0$ in $\mathbb{R}^d$, then $0\in \sigma(H)$.
\v{S}nol's theorem was generalized by many authors including
B.~Simon, see for example \cite{Cycon,S82} and \cite{Shubin}.

In \cite{kuchy1,kuchy2} the structure of the space of all
polynomially growing solutions of a periodic elliptic operator (or
a system) of order $m$ on an abelian cover of a compact Riemannian
manifold was studied. An important particular case of the general
results in \cite{kuchy1,kuchy2} is a real, second-order
$\mathbb{Z}^d$-periodic elliptic operator $P$ of the form
(\ref{P}) which is defined on $\mathbb{R}^d$. In this case, we can
use the information about positive solutions of such equations
described in Section~\ref{secperiod} and the results of
\cite{kuchy1} to obtain the precise structure and dimension of the
space of polynomially growing solutions.

\begin{definition} \label{defFermi} {\em 1. Let $N\geq 0$. We say that
{\em the Liouville theorem of order $N$ for the equation $Pu=0$}
holds true in $\mathbb{R}^d$, if the space $\mathrm{V}_N (P)$ of
solutions of the equation $Pu=0$ in $\mathbb{R}^d$ that satisfy
$\abs{u(x)}\leq C(1+\abs{x})^N$ for all $x\in \mathbb{R}^d$ is of
finite dimension.

2. The {\em Fermi surface} $F_P$ of the operator $P$ consists of
all vectors $\zeta\in \mathbb{C}^{d}$ such that the equation
$Pu=0$ has a nonzero {\em Bloch solution} of the form
$u(x)=e^{\mathrm{i}\zeta\cdot x}p(x)$, where
 $p$ is a $\mathbb{Z}^d$-periodic function.}
\end{definition}
For a general $\mathbb{Z}^d$-periodic elliptic operator $P$ of any
order, we have:
\begin{theorem}[\cite{kuchy1}]\label{thmliouv}
\begin{enumerate}
\item If the Liouville theorem of an order $N\geq 0$ for the equation $Pu=0$ holds true,
then it holds for any order.

\item The Liouville theorem holds true if and only if the number of points in the real Fermi
surface $F_{P}\cap \mathbb{R}^d$ is finite.
\end{enumerate}
\end{theorem}
For second-order operators with real coefficients, we have:
\begin{theorem}[\cite{kuchy1}]\label{T:non-self}Let $P$ be a $\mathbb{Z}^d$-periodic operator
on $\mathbb{R}^d$ of the form (\ref{P}) such that $\Lambda_{0}
\geq 0$. Then
\begin{enumerate}

\item The Liouville theorem holds vacuously if $\Lambda(0)> 0$, i.e., the equation
$Lu=0$ does not admit any nontrivial polynomially growing
solution.

\item If $\Lambda(0)= 0$ and $\Lambda_{0}>0$, then the Liouville theorem holds for
$P$, and
$$
\dim \mathrm{V}_N
(P)=\left(\begin{array}{c}d+N-1\\N\end{array}\right).
$$
\item If $\Lambda(0)= 0$ and $\Lambda_{0}=0$, then the Liouville theorem holds for
$P$, and
$$
\dim \mathrm{V}_N (P)=
\left(\begin{array}{c}d+N\\N\end{array}\right)-\left(\begin{array}{c}d+N-2\\N-2\end{array}\right),
$$
which is the dimension of the space of all harmonic polynomials of
degree at most $N$ in $d$ variables.

\item Any solution $u\in \mathrm{V}_N(P)$ of the equation
$Pu=0$ can be represented as
$$
u(x)=\sum\limits_{\abs{j}\leq N} x^j p_j(x)
$$
with $\mathbb{Z}^d$-periodic functions $p_j$.
\end{enumerate}
\end{theorem}

\mysection{Criticality theory for the $p$-Laplacian with potential
term}\label{secplap} Positivity properties of quasilinear elliptic
equations defined on a domain $\Omega\subset \mathbb{R}^d$, and in
particular, those with the $p$-Laplacian term in the principal
part, have been extensively studied over the recent decades (see
for example \cite{AH1,AH2,DKN,GS,HKM,V} and the references
therein).

Let $p\in(1,\infty)$, and let $\Omega$ be a general domains in
$\mathbb{R}^d$. Denote by $\Delta_p(u):=\nabla\cdot(|\nabla
u|^{p-2}\nabla u)$  the $p$-Laplacian operator, and let $V\in
L_\mathrm{loc}^\infty(\Omega)$ be a given (real) potential.
Throughout this section we always assume that \be \label{Q}
Q(u):=\int_\Omega \left(|\nabla u|^p+V|u|^p\right)\dx\geq 0\qquad
\forall u\in \core,\ee
that is, the functional $Q$ is nonnegative on $\core$. In
\cite{PT3}, K.~Tintarev and the author studied (sub)criticality
properties for positive weak solutions of the corresponding
Euler-Lagrange equation \be \label{groundstate}
 \frac{1}{p}Q^\prime
(v):=-\Delta_p(v)+V|v|^{p-2}v=0\quad \mbox{in }  \Omega,\ee  along
the lines of criticality theory for second-order linear elliptic
operators that was discussed in
sections~\ref{secpreliminaries}--\ref{secindef}.
\begin{definition}{\em We say that the functional $Q$ is {\em subcritical}
in $\Omega$ (or $Q$ is {\em strictly positive} in $\Omega$) if
there is a strictly positive continuous function $W$ in $\Omega$
such that \be\label{wsg}  Q(u)\geq  \int_\Omega W|u|^p\dx \qquad
\forall u\in\core.\ee }\end{definition}
\begin{definition}{\em We say that a sequence
$\{u_n\}\subset\core$  is a {\em null sequence}, if $u_n\geq 0$
for all $n\in\N$, and there exists an open set $B\Subset\Omega$
 such that
$\int_B|u_n|^p\dx=1$, and \be
\lim_{n\to\infty}Q(u_n)=\lim_{n\to\infty}\int_\Omega (|\nabla
u_n|^p+V|u_n|^p)\dx=0.\ee
 We say that a positive function $\varphi\in
C^1_{\mathrm{loc}}(\Omega)$ is a {\em ground state} of the
functional $Q$ in $\Omega$ if $\varphi$ is an
$L^p_{\mathrm{loc}}(\Omega)$ limit of a null sequence. If $Q\geq
0$, and $Q$ admits a ground state in $\Omega$, we say that the
functional $Q$ is {\em critical} in $\Omega$. The functional $Q$
is {\em supercritical} in $\Omega$ if $Q\ngeq 0$ on $\core$.}
\end{definition}
The following is a generalization of the Allegretto-Piepenbrink
theorem.
\begin{theorem}[see \cite{PT3}]\label{pos} Let Q be a functional of the form (\ref{Q}).
Then the following assertions are equivalent

(i) The functional $Q$ is nonnegative on $C_0^\infty(\Omega)$.

(ii) Equation (\ref{groundstate}) admits a global positive
solution.

(iii) Equation (\ref{groundstate}) admits a global positive
supersolution.
\end{theorem}
The definition of positive solutions of minimal growth in a
neighborhood of infinity in $\Omega$ in the linear case
(Definition~\ref{defminimalg}) is naturally extended to solutions
of the equation $Q'(u)=0$.
\begin{definition} \label{defminimalp}{\em
A positive solution $u$ of the equation $Q'(u)=0$ in $\Omega^*_j$
is said to be a {\em positive solution of the equation $Q'(u)=0$
of minimal growth in a neighborhood of infinity in} $\Omega$ if
for any $v\in C(\Omega_l^*\cup \partial  \Omega _l)$ with $l> j$,
which is a positive solution of the equation $Q'(u)=0$ in $\Omega
_l^*$, the inequality $u\le v$ on $\partial  \Omega _l$, implies
that  $u\le v$ on $\Omega _l^*$.}
\end{definition}

If  $1<p\leq d$, then for each $x_0\in \Omega$, any positive
solution $v$ of the equation $Q'(u)=0$ in a punctured neighborhood
of $x_0$ has either a removable singularity at $x_0$, or
 \begin{equation}\label{nonremovasymp}
  v(x)\sim\begin{cases}
    \abs{x-x_0}^{\alpha(d,p)} & p<d, \\
     -\log \abs{x-x_0} & p=d,
  \end{cases} \qquad \mbox{ as } x\to x_0,
\end{equation}
where $\alpha(d,p):=(p-d)/(p-1)$, and $f\sim g$ means that $
\lim_{x\to x_0}[{f(x)}/{g(x)}]= C$ for some $C>0$ (see \cite{GiS}
for $p=2$, and \cite{Serrin1,Serrin2,V,PT3} for $1<p\leq d$).
The following result is an extension to the $p$-Laplacian of
Theorem~\ref{thmmingr2}.
\begin{theorem}[\cite{PT3}]\label{thmmingr}
Suppose that $1<p\leq d$, and $Q$ is nonnegative on $\core$.
Then for any $x_0\in \Omega$ the equation $Q'(u)=0$ has (up to a
multiple constant) a unique positive solution $v$ in
$\Omega\setminus\{x_0\}$ of minimal growth in a neighborhood of
infinity in $\Omega$. Moreover, $v$ is either a global minimal
solution of the equation $Q'(u)=0$ in $\Omega$, or $v$ has a
nonremovable singularity at $x_0$.
\end{theorem}
The main result of this section is as follows.
\begin{theorem}[\cite{PT3}] \label{mainky3}
 Let $\Omega\subseteq\mathbb{R}^d$ be a domain,
$V\in L_\mathrm{loc}^\infty(\Omega)$, and $p\in(1,\infty)$.
Suppose that the functional $Q$ is nonnegative on $\core$. Then
\begin{itemize}
\item[(a)] The
functional $Q$ is either subcritical or critical in $\Omega$.

\item[(b)] If the functional $Q$ admits a ground state $v$, then $v$ satisfies
(\ref{groundstate}).

\item[(c)]
The functional $Q$ is critical in $\Omega$ if and only if
(\ref{groundstate}) admits a unique positive supersolution.

\item[(d)] Suppose that $1<p\leq d$. Then the functional $Q$
is critical (resp. subcritical) in $\Omega$ if and only if there
is a unique (up to a multiplicative constant) positive solution
$\varphi_0$ (resp. $G_Q^\Omega(\cdot,x_0)$) of the equation
$Q'(u)=0$ in $\Omega\setminus\{x_0\}$ which has minimal growth in
a neighborhood of infinity in $\Omega$ and has a removable (resp.
nonremovable) singularity at $x_0$.

\item[(e)]  Suppose that $Q$ has a ground state $\varphi_0$. Then there exists
a positive continuous function $W$ in $\Omega$, such that for
every $\psi\in C_0^\infty(\Omega)$ satisfying $\int_\Omega \psi
\varphi_0 \,\mathrm{d}x \neq 0$ there exists a constant $C>0$ such
that the following Poincar\'e type inequality holds:
 \be\label{Poinc}
 Q(u)+C\left|\int_{\Omega} \psi
 u\,\mathrm{d}x\right|^p \geq  C^{-1}\int_{\Omega} W|u|^p\,\mathrm{d}x \qquad
 \forall u\in C_0^\infty(\Omega).\ee
\end{itemize}
\end{theorem}
\begin{remarks}\label{remplapl2}{\em
1. Theorem \ref{mainky3} extends \cite[Theorem~1.5]{PT2} that
deals with the linear case $p=2$. The proof of Theorem
\ref{mainky3} relies on the {\em (generalized) Picone identity}
\cite{AH1,AH2}.

2. We call $G_Q^\Omega(\cdot,x_0)$ (after an appropriate
normalization) \emph{the positive minimal $p$-Green function of
the functional $Q$ in $\Omega$ with a pole at $x_0$}.

3. Suppose that $p=2$, and that there exists a function $\psi\in
L^2(\Omega)$ and $C\in \mathbb{R}$ such that
 \be\label{Poinc1}
     Q(u)+C\left|\int_{\Omega} \psi u\,\mathrm{d}x\right|^2\geq 0 \qquad
     \forall u\in C_0^\infty(\Omega),
     \ee
then the negative $L^2$-spectrum of $Q'$ is either empty or
consists of a single simple eigenvalue.
 }\end{remarks}

We state now several positivity properties of the functional $Q$
in parallel to the criticality theory presented in
sections~\ref{secpreliminaries}--\ref{secindef}.
For $V\in L^\infty_{\mathrm{loc}}(\Omega)$, we use the notation
\be Q_V(u):=\int_\Omega(|\nabla u|^p+V|u|^p)\dx\ee to emphasize
the dependence of $Q$ on the potential $V$.
\begin{proposition}
\label{monPot} Let  $V_j\in L^\infty_{\mathrm{loc}}(\Omega)$,
$j=1,2$. If $V_2\gneqq V_1$ and $Q_{V_1}\geq 0$ in $\Omega$, then
$Q_{V_2}$ is subcritical in $\Omega$.
\end{proposition}
\begin{proposition}
\label{monDom}
 Let $\Omega_1\subset\Omega_2$ be domains in $\R^d$ such that $\Omega_2\setminus\overline{\Omega_1}\neq\emptyset$.
Let $Q_V$ be defined on $C_0^\infty(\Omega_2)$.

1. If $Q_V\geq 0$ on $C_0^\infty(\Omega_2)$, then $Q_V$ is
subcritical in $\Omega_1$.

2. If $Q_V$ is critical in $\Omega_1$, then $Q_V$ is supercritical
in $\Omega_2$.
\end{proposition}
\begin{proposition}\label{Prop2}
 Let $V_0, V_1\in L^\infty_{\mathrm{loc}}(\Omega)$,
$V_0\neq V_1$. For $s\in \mathbb{R}$ we denote \be
Q_s(u):=sQ_{V_1}(u)+(1-s)Q_{V_0}(u),\ee and suppose that
$Q_{V_j}\geq 0$ on $\core$ for $j=0,1$.

Then the functional $Q_s\geq 0$ on $\core$ for all $s\in[0,1]$.
Moreover, if $V_0\neq V_1$, then $Q_s$ is subcritical in $\Omega$
for all $s\in(0,1)$.
\end{proposition}
\begin{proposition}\label{strictpos}
Let $Q_V$ be a subcritical in $\Omega$. Consider $V_0\in
L^\infty(\Omega)$ such that $V_0\ngeq 0$ and $\supp
V_0\Subset\Omega$. Then there exist $0<\tau_+<\infty$, and
$-\infty\leq \tau_-<0$ such that $Q_{V+sV_0}$ is subcritical in
$\Omega$ for $s\in(\tau_-,\tau_+)$, and $Q_{V+\tau_+ V_0}$ is
critical in $\Omega$. Moreover, $\tau_-=-\infty$ if and only if
$V_0\leq 0$.
\end{proposition}
\begin{proposition}\label{propintcond}
Let $Q_V$ be a critical functional in $\Omega$, and let
$\varphi_0$ be the corresponding ground state. Consider $V_0\in
L^\infty(\Omega)$ such that $\supp V_0\Subset\Omega$. Then there
exists $0<\tau_+\leq\infty$ such that $Q_{V+sV_0}$ is subcritical
in $\Omega$ for $s\in(0,\tau_+)$ if and only if
 \be\label{intcond}
  \int_\Omega
V_0(x)\varphi_0(x)^p\dx>0.
 \ee
\end{proposition}
\begin{center}
{\bf Acknowledgments} \end{center} The author expresses his
gratitude to F.~Gesztesy and M.~Murata for their valuable remarks.
The author is also grateful to the anonymous referee for his
careful reading and useful comments. This work was partially
supported by the Israel Science Foundation founded by the Israeli
Academy of Sciences and Humanities, and the Fund for the Promotion
of Research at the Technion.

\end{document}